\documentclass[12pt]{article}
\usepackage[margin=1in]{geometry}
\usepackage{amsmath,amssymb,amsthm}
\usepackage{caption}
\usepackage{mathtools}
\usepackage{booktabs}
\usepackage{setspace}
\usepackage{natbib}
\usepackage{titlesec}
\usepackage{rotating}
\usepackage{float}
\usepackage[colorlinks=true, allcolors=blue]{hyperref}
\usepackage[toc,page]{appendix}
\usepackage{appendix}
\usepackage{xcolor}
\usepackage[normalem]{ulem} 

\newtheorem{theorem}{Theorem}
\newtheorem{lemma}{Lemma}
\newtheorem{proposition}{Proposition}
\newtheorem{corollary}{Corollary}
\newtheorem{assumption}{Assumption}
\newtheorem{hypothesis}{Hypothesis}
\theoremstyle{remark}
\newtheorem{remark}{Remark}

\newcommand{\E}{\mathbb{E}}
\newcommand{\Prob}{\mathbb{P}}

\newcommand{\Cov}{\operatorname{Cov}}
\newcommand{\plim}{\overset{p}{\longrightarrow}}
\newcommand{\dlim}{\overset{d}{\longrightarrow}}
\newcommand{\Lam}{\Lambda}
\newcommand{\F}{\mathcal{F}}
\newcommand{\Gsig}{\mathcal{G}}
\newcommand{\QC}{\widehat{QC}}        
\newcommand{\QCd}{\widehat{QC}^{\,0}} 
\newcommand{\Rt}{\widetilde{R}}
\newcommand{\rt}{\widetilde{r}}

\title{\vspace*{-2cm} Inference for Realized Quadratic Covariation of $R$ and $R^2$ as a Measure of Return Asymmetry}

\author{Haozhe Jiang,$^a$ Ostap Okhrin,$^{b}$ Michael Rockinger,$^{c,*}$}

\date{August 2026}
\begin{document}
\maketitle
\begin{abstract}

We propose a nonparametric test of return asymmetry based on the realized quadratic covariation between a cumulative return path and its square using high-frequency data. Our null hypothesis is that the volatility path is independent of the Brownian motion driving returns. Under this null and a martingale condition, the infill limit of the statistic is conditionally Gaussian given the volatility path, with a conditional variance determined by integrated variance. Studentization by a feasible variance estimator yields an asymptotically standard normal statistic without specifying a parametric model for volatility. 
Simulations under Heston and Bates dynamics show near-nominal per-window size, quantify power gains from independent-window aggregation, and assess sensitivity to jump removal and volatility-dependent drift.
An application to 165 U.S. stocks over 1998--2025 finds that 
rejections concentrate in negative asymmetry, consistent with the leverage effect.


\vspace*{0.07cm}

\vspace*{0.03cm}
\end{abstract}
\clearpage

\vspace*{0.5cm}
\noindent$\mbox{}^{a}\,${\footnotesize 
Technische Universität Dresden. ``Friedrich List'' Faculty of Transportation. Chair of Applied Statistics. E-mail: \url{haozhe.jiang@tu-dresden.de}.
}

\noindent$\mbox{}^{b}\,${\footnotesize 
Technische Universität Dresden. ``Friedrich List'' Faculty of Transportation. Chair of Applied Statistics.  \\ Center for Scalable Data Analytics and Artificial Intelligence (ScaDS.AI) Dresden/Leipzig, Germany. E-mail: \url{ostap.okhrin@tu-dresden.de}.
}

\noindent$\mbox{}^{c}\,${\footnotesize 
University of Lausanne. Faculty of
Business and Economics, CH 1015 Lausanne, Switzerland. E-mail: \url{michael.rockinger@unil.ch}.}

\vspace*{0.5cm}
{\footnotesize  \noindent * Corresponding author. The views expressed in this paper are those of the authors and are not necessarily reflective of views from their employers. We report no conflict of interest. We are most grateful to the Deutsche Forschungsgemeinschaft (DFG, German Research Foundation) who partially funded the project under number 497262150.}

\clearpage\newpage
\section{Introduction}

The \emph{leverage effect}, the correlation between returns and volatility, is a common term used to evaluate return asymmetry \citep{EngleNg1993, Zakoian1994, WuXiao2002}. Building on this notion, \citet{WangMykland2014} and \citet{AitSahalia2017} estimate the leverage effect by deriving the asymptotic distribution of the realized quadratic covariation between the return and (a function of) the volatility process.
In this paper, with a focus on return asymmetry, we instead study the realized quadratic covariation between the cumulative returns $R$ and its running square~$R^2$, denoted
$\QC_T$. As a sum of products of increments, $\QC_T$ converges in probability
to the quadratic covariation $[R,R^2]_T$ over any sequence of partitions with
mesh tending to zero, requiring nothing beyond the semimartingale property
\citep{jacod1987limit,Barndorff-N&Shephard2004a}. The probability limit
decomposes into an integral against the quadratic variation and a cubed-jump
sum, $[R,R^2]_T = 2\int_0^T R_{u-}\,d[R]_u + \sum_{u \le T}(\Delta R_u)^3$,
both of which carry information about the asymmetry of the return path. The first term captures the dependence between the return path and its volatility innovation, and up to the constant $1.5$ it coincides with the realized third cumulant of a low-frequency window return for a continuous semimartingale \citep{NP2020, choelee}. The second term describes the jump component in the realized skewness \citep{Amaya2016}. 
Thus, to develop a test of return asymmetry, we propose to use the realized quadratic covariation based on high-frequency intraday data of the stock market and obtain a test-statistic, whose limit distribution is ensured by the general semi-martingale theory. The only remaining question is what distribution theory is available for inference.

Specializing to the continuous stochastic volatility semimartingale class
($\mathcal{SVSM}^c$), \citet{Barndorff-N&Shephard2004a,Barndorff-N&Shephard2004c} provide the distribution theory of the realized
covariation matrix. Applied to the present pair $(R,R^2)$, their central limit
theorem gives
\begin{equation}\label{eq:bns_intro}
\frac{\QC_T - L}{\sqrt{H_n}} \;\dlim\; \mathcal{N}(0,1),
\end{equation}
where $L = 2\int_0^T R_u\,d\Lambda_u$ is the \emph{pathwise} limit of $\QC_T$ and
$H_n$ is a feasible estimator of the variance of the estimation error
$\QC_T - L$. Their theorem demonstrates that, standardized by the infill
variance, the discrepancy between the realized statistic and its pathwise
limit is asymptotically standard normal. So their objective is to quantifying the statistical deviation between the realized quadratic covariance from a pathwise generated benchmark. The same pathwise logic underlies the high-frequency literature on the leverage effect. \citet{WangMykland2014}, \citet{AitSahalia2017}, and \citet{KalninaXiu2017} estimate the quadratic covariation (or correlation) between returns and volatility and derive infill central limit theorems for the estimation error around its realized value. 
However, our purpose is different: \emph{we aim
to test whether the mechanism generating returns is asymmetric.} 
For the quadratic-covariation functional considered here, we construct a statistic of the form 
\begin{equation}\label{eq:mech_intro}
\frac{\QC_T - \E[L]}{\sqrt{\widehat{V}_n}},
\end{equation}
where $\widehat V_n$ estimates the conditional variance of $L$ given the volatility path. This variance is of order $O(1)$, whereas the infill error variance estimated by $H_n$ is of order $O_p(\Delta_n^{1/2})$.
%
Under the volatility--Brownian independence null and the martingale condition, $L$ has conditional mean zero. The statistic therefore reduces to $\widehat{QC}_T/\sqrt{\widehat V_n}$.
Since both the numerator
and the denominator have changed, a new feasible variance estimator and a
corresponding asymptotic normality result are needed, which is what the theory
of this paper develops.

We develop that test. Jumps are first removed with the nonparametric method of
\citet{LeeMyklandJump}, so that inference proceeds within the $\mathcal{SVSM}^c$ class. After filtering, the mean of the statistic is nonzero only if
return--variance dependence is present in the continuous part. The null hypothesis is volatility--Brownian independence, i.e., the volatility path is independent of the Brownian motion driving returns, which is a sufficient condition for the absence of return asymmetry. 
Conditionally on the volatility path, the limit of the realized quadratic
covariation is then exactly Gaussian with a variance that is a functional of
integrated variance alone under the martingale assumption. Studentizing by a feasible estimator of that conditional variance yields a statistic that is asymptotically standard normal under the null, regardless of the level, persistence, or memory of volatility, because the volatility path is \emph{divided out rather than modeled}. We then relax the martingale assumption to a window-constant drift and show that asymptotic normality is preserved once a demeaned version of the test statistic is used with a bridge correction. 
A single window of high-frequency observations, however, carries only one realization of the limit but little information about its mean. This difficulty echoes the \emph{leverage effect puzzle} documented by \citet{AitSahaliaFanLi2013}, who show how hard the leverage correlation is to recover nonparametrically at high frequency.
In this paper, power is obtained by aggregating the studentized statistics over $K$ non-overlapping windows under the assumption of stationary and ergodic window blocks. 
A Monte Carlo study under the \citet{heston1993closed} and \citet{bates1996jumps} stochastic differential equations validates the theory. To our knowledge, this is the first procedure that uses the realized quadratic covariation $[R,R^2]$, the high-frequency counterpart of the realized third cumulant, for inference on the return-generating mechanism rather than on the realized path. Its size is correct regardless of the volatility dynamics, and it makes explicit how much aggregation is required for power. The procedure complements existing nonparametric leverage estimators, as \citet{AitSahaliaFanLi2013} and \citet{KalninaXiu2017}, by providing inference for a third-order return-asymmetry functional without specifying parametric volatility dynamics.

The paper proceeds as follows. 
Section~\ref{sec:background} introduces the semimartingale framework and the realized quadratic covariation, and distinguishes inference on its realized limit from inference under the structural independence null.
Section~\ref{sec:test} then formalizes our null as Hypothesis~\ref{hyp:null}, develops the
conditional-Gaussian inference theory under the martingale condition, and proposes the aggregate mean test for power.
Section~\ref{sec:empimpl} relaxes the martingale assumption by a window-constant drift, and constructs the operational demeaned statistic with a bridge-corrected variance estimator. Section~\ref{sec:mc} reports the simulation
study, jump-free and jump-present, including the closed-form Heston
benchmark, the aggregation-based power analysis, and the jump-robustness
experiments. In Section~\ref{sec:realworld}, we apply the proposed test to the real-world financial market. Section~\ref{sec:conclusion} concludes. Proofs and derivations are collected in
the Appendix.

\section{Model and Background}\label{sec:background}

\subsection{Setup: process, sampling, and the realized quadratic
covariation}\label{sec:setup}

Fix a window $[0,T]$ (in the empirical implementation, one month of trading
time) and let the log-price $X$ of an asset be a semimartingale ($\mathcal{SM}$) on a
filtered probability space $(\Omega,\F,(\F_t)_{t\ge 0},\Prob)$. Following
\citet{Barndorff-N&Shephard2004a, Barndorff-N&Shephard2004c}, recall that
every $X\in\mathcal{SM}$ decomposes as
\begin{equation}\label{eq:smdecomp}
X \;=\; X_0+A+M ,
\end{equation}
with $A$ of finite variation and $M$ a local martingale. 
We work at two levels of generality: the object of study, i.e.,
quadratic covariation, is defined and consistently estimated with no
structure beyond $X\in\mathcal{SM}$ (Section~\ref{sec:consistency}),
whereas the distribution theory requires the continuous It\^o
specialization of Section~\ref{sec:svsm}.

The window is sampled at $n$ equispaced times $t_j = j\Delta_n$,
$j=0,1,\dots,n$, with $\Delta_n = T/n$. All limits below are taken as
$n\to\infty$ with $T$ fixed. Intraday returns are
$r_j = X_{t_j}-X_{t_{j-1}}$, and we write $R_t=X_t-X_0$ for the cumulative
return path. 
Throughout we replace the time index along a path or a filtration by its index.
Hence
$Z_j:=Z_{t_j}$, so that $R_j=X_j-X_0$ and $\F_j=\F_{t_j}$. Increments over
the $j$-th interval $(t_{j-1},t_j]$ are written $\Delta Z_j=Z_j-Z_{j-1}$,
so that $\Delta R_j=R_j-R_{j-1}=r_j$ and $R_{j-1}$ is the level at which that interval
begins. Time instances keep the symbol $t_j$.

The object of study is the realized quadratic covariation between the
return path and its square,
\begin{equation}\label{eq:QCraw}
\QC_T \;=\; \sum_{j=1}^{n}\Delta R_j\,\Delta\bigl(R^2\bigr)_j
\;=\;\sum_{j=1}^{n}\Bigl\{\,2\,R_{j-1}\bigl(\Delta R_j\bigr)^2
+\bigl(\Delta R_j\bigr)^3\Bigr\},
\end{equation}
the second equality following from
$R_j^2-R_{j-1}^2=2R_{j-1}\Delta R_j+(\Delta R_j)^2$, which is
where the level $R_{j-1}$ enters. Equivalently, $\QC_T$ is the off-diagonal entry of the realized covariation matrix of the bivariate process $(R,R^2)$, the object to which the limit
theory of Section~\ref{sec:bns} applies.
Up to the constant
$3/2$, $\QC_T$ is the quadratic-covariation component of the realized third
cumulant of the window return \citep{NP2020}. Since the hypotheses tested below set a mean
to zero, the constant is irrelevant. Section~\ref{sec:demean} introduces
the \emph{demeaned} version $\QCd_T$, computed from within-window demeaned
returns, which is the operational statistic of the empirical work.

\subsection{Quadratic covariation in the general semimartingale
framework}\label{sec:consistency}

The probability limit of $\QC_T$ is available in full generality, by the
theory of quadratic variation and with no assumption beyond
$X\in\mathcal{SM}$, in particular also with jumps present.

\begin{lemma}[Consistency of the realized quadratic covariation;
\citet{jacod1987limit,Barndorff-N&Shephard2004a}]\label{lem:limitraw}
Let $R\in\mathcal{SM}$. Throughout, $[R]$ denotes the \emph{quadratic
variation} of $R$: the non-decreasing c\`adl\`ag process to which the realized
variance along the grid converges in probability,
\[
\sum_{t_j\le t}\bigl(\Delta R_j\bigr)^2\;\plim\;[R]_t\qquad(n\to\infty),
\]
which exists for every $\mathcal{SM}$
\citep[Theorem~I.4.47]{jacod1987limit}. The quadratic covariation
$[R,R^2]$ of $R$ and $R^2$ is defined in the same way, with the summand
$\bigl(\Delta R_j\bigr)^2$ replaced by $\Delta R_j\,\Delta\bigl(R^2\bigr)_j$.
Then, as $n\to\infty$,
\[
\QC_T \;\plim\; \bigl[R,R^2\bigr]_T
\;=\; 2\int_0^T R_{u-}\,d[R]_u
\;+\;\sum_{0<u\le T}\bigl(\Delta R_u\bigr)^{3} .
\]
If moreover $R$ is continuous with $[R]=\Lam$, the jump sum vanishes and
\[
\QC_T \;\plim\; L \;:=\; 2\int_0^T R_u\, d\Lam_u ,
\]
where the dependence of $L$ on the fixed window length $T$ is suppressed
in the notation.
\end{lemma}

\begin{proof}
For any pair of semimartingales, the sum of products of increments
converges in probability to the quadratic covariation over any sequence of
partitions with mesh tending to zero \citep[Theorem~I.4.47 and I.4.52]{jacod1987limit}. This is the standard convergence of realized covariation reviewed in \cite{Barndorff-N&Shephard2004a}, and it holds for every semimartingale, with or without
jumps, and in particular without any condition on the dependence between
the volatility and the driving Brownian motion. Applying it to the pair
$(R,R^2)$ gives the first limit once $[R,R^2]_T$ is identified: It\^o's
formula for semimartingales yields
$d(R^2)=2R_{-}\,dR+d[R]$, whence
$[R,R^2]_T=2\int_0^T R_{u-}\,d[R]_u+\sum_{u\le T}(\Delta R_u)^3$, the jump
sum arising from the covariation of the jump parts. Under continuity the
sum is empty, $R_{u-}=R_u$ and $[R]=\Lam$, which is the second display.
\end{proof}

Consistency is thus the one layer of the theory that requires no assumptions beyond the semimartingale property. What it does not deliver is a distribution: neither the sampling error of $\QC_T$
around $[R,R^2]_T$ nor the law of the limit itself is available in
$\mathcal{SM}$.

\subsection{Specialization: continuous stochastic volatility
semimartingales and jump removal}\label{sec:svsm}

A central limit theory is not available in the general semimartingale
class, and \cite{Barndorff-N&Shephard2004c} accordingly specialize. We
follow them and assume that the log-price $X$ is an It\^o semimartingale,
\begin{equation}\label{eq:sm}
X_t \;=\; X_0 \;+\; \int_0^t \alpha_u\,du \;+\; \int_0^t \sigma_u\,dW_u
\;+\; J_t ,
\end{equation}
where $W$ is a standard Brownian motion, $\alpha$ is a c\`adl\`ag drift
with $\int_0^T\alpha_u^2\,du<\infty$, $\sigma$ is a c\`adl\`ag adapted
spot-volatility process, and $J$ is a finite-activity jump component. This
is a strict subclass of $\mathcal{SM}$: the drift and the jump compensator
are required to admit densities with respect to Lebesgue time, and the
continuous martingale part to be a stochastic integral against Brownian
motion, so that $[X]^c_t=\int_0^t\sigma_u^2\,du$. Setting $J\equiv0$ in
\eqref{eq:sm} gives the class which
\cite{Barndorff-N&Shephard2004c} named as a continuous
stochastic volatility semimartingale, $\mathcal{SVSM}^c$, which is the
class in which the limit theory of Section~\ref{sec:bns} operates. It is worth being explicit that the
continuity restriction applies to the \emph{price}: the volatility process
is required only to be c\`adl\`ag, so $\sigma$ may itself jump.

Empirically, the jump component is removed rather than assumed away.
Returns flagged as jumps by a high-frequency filter are set to zero. The conditional-Gaussian theory is developed for the continuous component and applies to filtered returns when filtering preserves the relevant infill limits; Remark~\ref{rem:jumps} discusses classification error.
Under finite-activity jumps, such conditions are satisfied almost surely by consistent threshold procedures, including the method of \citet{mancini2009}. In the simulation and real-world study, we use \citet{LeeMyklandJump} as the baseline implementation. At finite sampling frequencies, neither procedure recovers the continuous increment on a jump-containing interval exactly, because zeroing that observation also discards its diffusive component. However, the cumulative discarded diffusive contribution is asymptotically negligible. 
This step is what
makes the specialization applicable to equity data, and it is also what
gives the statistic its economic content, as the interpretation below
explains.

Write $\Lam_t=\int_0^t\sigma_u^2\,du$ for the integrated variance, the
quadratic variation of the continuous martingale part of $X$, and let
\[
\Gsig \;=\; \sigma\bigl((\sigma_u)_{u\in[0,T]}\bigr)
\]
denote the $\sigma$-field generated by the volatility path over the window.
Every statement below made ``conditionally on the volatility path'' is a
statement conditional on $\Gsig$. With jumps removed, the second part of
Lemma~\ref{lem:limitraw} applies and $\QC_T\plim L=2\int_0^TR_u\,d\Lam_u$.

\begin{assumption}[Continuous stochastic 
volatility structure]\label{ass:sv}
On $[0,T]$, $X$ is of the form \eqref{eq:sm} with $J\equiv0$, and $\sigma$
is c\`adl\`ag with
$\underline{\sigma}\le\sigma_u\le\overline{\sigma}$ for constants
$0<\underline{\sigma}\le\overline{\sigma}<\infty$.
\end{assumption}

\begin{assumption}[Martingale condition]\label{ass:mart}
$\alpha\equiv 0$ on $[0,T]$, so that $R$ is a continuous local martingale
with $[R]=\Lam$.
\end{assumption}

Assumption~\ref{ass:mart} is imposed throughout the inference theory of
Section~\ref{sec:mart}. It is restrictive for asset prices, and
Section~\ref{sec:demean} shows that the demeaned statistic computed from
returns carrying any window-constant drift coincides exactly, pathwise,
with the statistic computed from the martingale part, so the operational
results are invariant to such drifts.

\paragraph{Interpretation: return--variance dependence and the leverage
effect.} Although $\QC_T$ is formally a covariation between $R$ and its
square, $\QC_T$ is not a simple dependence
measure between these two processes as obviously $R^2$ is a functional of $R$. 
Its economic content is
\emph{return--variance dependence}. Throughout, by \emph{return--variance dependence} we mean dependence between the
Brownian motion driving returns and the variance path, as opposed to the
trivial dependence of the return path on the variance path through its
scale.
The expansion in \eqref{eq:QCraw} shows
that the dominant summand, $2R_{j-1}(\Delta R_j)^2$, pairs the past
cumulative return with the current squared return, so that the statistic
aggregates the covariance between return levels and subsequent realized
variance. 
Under Assumptions~\ref{ass:sv}--\ref{ass:mart}, the results below imply
\begin{equation}\label{eq:meanQC}
\lim_{n\to\infty}\E\bigl[\QC_T\bigr] \;=\;
2\int_0^T \Cov\bigl(R_u,\,\sigma_u^2\bigr)\,du ,
\end{equation}
the integrated covariance between the cumulative return and the spot
variance.\footnote{This combines the probability limit $\QC_T\plim L$ of
Lemma~\ref{lem:limitraw} with convergence of the means: the boundedness of
$\sigma$ in Assumption~\ref{ass:sv} makes $\{\QC_T\}_n$ uniformly
integrable, so that
$\E[\QC_T]\to\E[L]=2\int_0^T\Cov(R_u,\sigma_u^2)\,du$.} 
%
A nonzero integrated return--variance covariance implies dependence between the volatility path and the Brownian motion driving returns, and its sign agrees with that of the third cumulant of the martingale return. In the Heston setting considered in the simulation study, this integrated covariance is proportional to the instantaneous correlation $\rho$, i.e., the \emph{leverage effect}. The resulting test is calibrated under volatility--Brownian independence and targets directional departures associated with return asymmetry.


\subsection{Background: the Barndorff-Nielsen \& Shephard limit
theory}\label{sec:bns}

Lemma~\ref{lem:limitraw} supplied the first layer of the
realized-covariation theory. Two further layers complete it, and both
require the $\mathcal{SVSM}^c$ specialization of Section~\ref{sec:svsm}.

\emph{(i) The central limit theorem for the pathwise estimation error.}
For the realized covariation of an $\mathcal{SVSM}^c$ pair with
square-integrable drift, \cite{Barndorff-N&Shephard2004a} establish a
central limit theorem for the difference between the realized covariation
and its pathwise limit, at rate $\Delta_n^{-1/2}$. Their original argument
conditions on the volatility path and therefore assumes no leverage, i.e.,
independence of the covolatility from the driving Brownian motion.
\cite{Barndorff-N&Shephard2004c}, based on \cite{JacodProtter1998}, prove that the same conclusion holds under arbitrary
diffusion-based leverage, so that no restriction relating the volatility to
$W$ is required. Under the regularity conditions of the cited central limit theorem, its specialization to the pair $(R,R^2)$ gives 
\begin{equation}\label{eq:bnsclt}
\Delta_n^{-1/2}\bigl(\QC_T - L\bigr) \;\dlim\; \mathrm{MN}(0,\Pi_T),
\end{equation}
where the convergence is stable in law and the limiting mixed normal is driven by a Brownian motion independent of $X$, the feature that makes feasible studentization legitimate \citep{Barndorff-N&Shephard2004c}.
It\^o's formula gives $d(R^2)=\sigma_u^2\,du+2R_u\sigma_u\,dW_u$,
so its drift and diffusion coefficients are $(0,\sigma_u^2)'$ and
$(\sigma_u,2R_u\sigma_u)'$, both c\`adl\`ag with
$\int_0^T\sigma_u^4\,du<\infty$ under Assumption~\ref{ass:sv} and~\ref{ass:mart}.
Specializing the Barndorff-Nielsen--Shephard asymptotic covariance to this
diffusion coefficient gives the conditional variance in closed form,
$\Pi_T \;=\; 8\int_0^T R_u^2\,\sigma_u^4\,du$.
The theory is completed by a feasible estimator of $\Pi_T$.
\cite{Barndorff-N&Shephard2004a} show how the asymptotic covariance
array of the realized covariation matrix can be estimated from
high-frequency data as follows.
\begin{equation}\label{eq:Hn}
H_n \;=\; \sum_{j=1}^{n} x_j^2 \;-\; \sum_{j=1}^{n-1} x_j x_{j+1},
\qquad
\Delta_n^{-1}H_n \;\plim\; \Pi_T ,
\end{equation}
formed from the summands $x_j=\Delta R_j\,\Delta(R^2)_j$ of $\QC_T$
itself and corrected for their first-order autocorrelation. Thus
$H_n=O_p(\Delta_n)$.

\emph{(ii) Asymptotic normality of the standardized statistic.}
Combining \eqref{eq:bnsclt}--\eqref{eq:Hn} with stability of the
convergence,
\begin{equation}\label{eq:bnst}
\frac{\QC_T - L}{\sqrt{H_n}} \;\dlim\; \mathcal{N}(0,1).
\end{equation}
The standardized pathwise estimation error is asymptotically standard
normal, and Gaussian confidence intervals for the realized quantity $L$ follow path by path.

\subsection{Pathwise and structural hypotheses}\label{sec:twonulls}


Together, Lemma~\ref{lem:limitraw} and Section~\ref{sec:bns} constitute a
\emph{measurement} theory: they identify what $\QC_T$ estimates, the realized path
quantity $L=2\int_0^T R_u\,d\Lam_u$, and how precisely it is measured.
Whether the standardized statistic \eqref{eq:bnst} can test a hypothesis
of the mechanism level is a separate question, and the answer turns on
distinguishing
\[
H_0^{\mathrm{path}}:\; L=0 \text{ for the given window},
\qquad\text{versus}\qquad
H_0^{\mathrm{mech}}:\; \text{$\sigma\perp W$ (implying $\E[L]=0$)},
\]
i.e., a statement about one realized path versus a statement about the
return-generating \emph{mechanism}. Under $H_0^{\mathrm{path}}$ we test whether the $L$ is zero for one particular window, whereas under $H_0^{\mathrm{mech}}$ the target is whether the Brownian motion driving returns is independent of the volatility path, and thus the limit $L$, a continuous random variable, has mean zero, i.e., $\E[L]=0$. 

Replacing the unobservable centering $L$ in \eqref{eq:bnst} by the
hypothesized value zero yields the $t$-statistic $\QC_T/\sqrt{H_n}$. Adding and
subtracting $L$ splits the statistic into the standardized estimation
error and an offset,
\begin{equation}\label{eq:divergence}
\frac{\QC_T}{\sqrt{H_n}}
\;=\;\underbrace{\frac{\QC_T-L}{\sqrt{H_n}}}_{\dlim\;\mathcal{N}(0,1)}
\;+\;\underbrace{\frac{L}{\sqrt{H_n}}}_{
=\;\Delta_n^{-1/2}\,L/\sqrt{\Pi_T}\;+\;o_p(\Delta_n^{-1/2})} .
\end{equation}
The first term is the object the limit theory delivers and is $O_p(1)$. In the
second term the numerator, as a property of the realized path, is $O_p(1)$,
while the denominator is $O_p(\Delta_n^{1/2})$ by \eqref{eq:Hn}. 
%
%
Under $H_0^{\mathrm{mech}}$, the conditional Gaussian law established in Section~\ref{sec:test} implies $\Prob(L=0)=0$. Consequently, $\lvert\widehat{QC}_T/\sqrt{H_n}\rvert$ diverges in probability, and the rejection probability tends to one under this null.
Equivalently, in confidence-interval form: \eqref{eq:bnst} gives
$L\in\QC_T\pm1.96\sqrt{H_n}$, an interval of width $O_p(\Delta_n^{1/2})$
collapsing onto a nonzero limit, so eventually it cannot contain the
origin.
%
This rejection concerns the realized value of $L$ and does not, by itself, provide evidence against volatility--Brownian independence.
In the simulations of Section~\ref{sec:mc} we will encounter nearly always rejection, which is exactly what our finding predicts.
%
We therefore normalize $\widehat{QC}_T$ by an estimator of the conditional standard deviation of $L$ given the volatility path, which remains of order $O(1)$ as sampling becomes finer.
That is the program of Section~\ref{sec:test}. 

The two theories are complements operating at
different stochastic orders: the classical theory of \citet{Barndorff-N&Shephard2004a}
answers \emph{what is the realized quadratic covariation in this window, and
how precisely is it measured?}, while the theory below answers \emph{given a window of observations, is the observed quadratic covariation compatible with volatility--Brownian independence?}.

\section{Contribution}\label{sec:test}

\subsection{The null hypothesis}\label{sec:null}

The theory of this section provides feasible inference for the mechanism
hypothesis left untested by the standardized statistic of
Section~\ref{sec:bns}. We first state the null formally.
\setcounter{hypothesis}{-1}
\begin{hypothesis}[Volatility--Brownian independence, $\sigma\perp W$]\label{hyp:null}
The volatility path $\sigma=(\sigma_u)_{u\in[0,T]}$ is independent of the
Brownian motion $W$ driving returns.
\end{hypothesis}

We refer to Hypothesis~\ref{hyp:null} as volatility--Brownian independence, or $H_0$, throughout. It is the null of this paper: within the continuous class $\mathcal{SVSM}^c$, it rules out any dependence between the Brownian motion driving returns and the volatility path, and together with the symmetry of Brownian increments it implies that, conditionally on the
volatility path, returns are exactly symmetric. Specialized to a leverage-diffusion model such as Heston, it reduces to $\rho=0$. 

The alternative $H_1$ is its complement, i.e., $H_1:\ \sigma\not\perp W$, so that the volatility path carries information about the Brownian motion driving returns.

The quadratic covariation $L$ does not, however, provide an omnibus measure of
arbitrary dependence between $W$ and $\sigma$. Under
Assumption~\ref{ass:mart}, its population location is
the signed integrated return--variance dependence
\begin{equation}\label{eq:integraCOV}
\E[L]
= 2\int_0^T \Cov(R_u,\sigma_u^2)\,du.
\end{equation}
Hence $H_0$ implies $\E[L]=0$, but the converse need not hold: dependence
between return and volatility innovations may be present while positive and
negative contributions to the covariance integral offset each other.
Accordingly, the statistic measures net directional third-order return asymmetry
rather than arbitrary return--volatility dependence. A nonzero value
$\E[L]\neq0$ therefore provides evidence against $H_0$, whereas $\E[L]=0$
does not by itself establish volatility--Brownian independence. In the Heston
model considered in Section~\ref{sec:hestonbench}, this distinction disappears because
$\E[L]$ is proportional to $\rho$ and vanishes if and only if $\rho=0$. Return--variance dependence constitutes the alternative, while, owing to the nature of quadratic covariation, the proposed test is designed to detect directional return--variance dependence that produces a nonzero signed integrated covariance.
 
Under $H_0$, conditionally on $\Gsig$ the return process $R$
is a continuous Gaussian martingale with deterministic variance function
$\Lam$. This conditional structure drives every distributional result
below. Unlike Assumptions~\ref{ass:sv}--\ref{ass:mart}, $H_0$ is a
hypothesis rather than a maintained condition: the consistency layer
(Lemmas~\ref{lem:limitraw} and~\ref{lem:limitdm}) holds with or without
it, the null distribution theory is derived under it, and power against
its violation is characterized in Corollary~\ref{cor:power} and quantified in
Section~\ref{sec:mc}.

\subsection{Inference under the martingale condition}\label{sec:mart}

The consistency layer, Lemma~\ref{lem:limitraw}, is classical. What must
be built is the distribution theory for the $O(1)$ limit $L$ itself,
which is not covered by the classical theory
(Section~\ref{sec:twonulls}) and is where $H_0$ enters. 
The proofs of all results in this section are provided in Appendix~\ref{app:proofs-mart}.

\begin{lemma}[Riemann approximation of integrals against $d\Lam$ under $H_0$]
\label{lem:riemann}
Let Assumptions~\ref{ass:sv} and $H_0$ hold and let
$(f_u)_{u\in[0,T]}$ be a continuous process, adapted to $\F_u\vee\Gsig$,
with $\E\bigl[\sup_{u\le T}f_u^2\mid\Gsig\bigr]<\infty$ a.s. Then
\[
\sum_{j=1}^n f_{j-1}\,(\Delta R_j)^2 \;\plim\; \int_0^T f_u\,d\Lam_u ,
\qquad n\to\infty .
\]
\end{lemma}

For integrands adapted to $(\F_u)$ the conclusion is the classical
convergence of weighted realized quadratic variation, valid for every
continuous semimartingale and without $H_0$
\citep[Theorem~I.4.47 and I.4.54]{jacod1987limit}. 

The content of the lemma is the
extension to weights that are measurable with respect to
$\F_u\vee\Gsig$, such as $\Lam_T-\Lam_u$ in
Proposition~\ref{prop:condgaussraw} and the bridge constant $c$ in
Theorem~\ref{thm:consistency}, which anticipate the volatility path
over the whole window and therefore lie outside the classical theory.
This extension is where $H_0$ enters: under independence of $\sigma$ and
$W$, enlarging the filtration by $\Gsig$ leaves $W_t$ still a Brownian motion and the martingale part ($M_t=\int_0^t \sigma_u\,dW_u$) of $R$ still a martingale, so that
conditionally on $\Gsig$ the anticipative weights are deterministic and
an elementary argument applies. Under leverage the enlargement destroys
the martingale property, i.e., $\Gsig$ carries information about the future increments of $W_t$,
$M_t$ acquires a drift in the enlarged filtration, and no such reduction is available.


Lemma~\ref{lem:riemann} evaluates the weight at the left endpoint
$t_{j-1}$ of each sampling interval, as is natural for It\^o-type sums.
The feasible variance estimators below, however, weight the $j$-th
squared return by quantities evaluated at the right endpoint $t_j$, such
as $\widehat\Lam_n(T)-\widehat\Lam_n(j)$ in \eqref{eq:Vhatraw} and $j/n$
in the bridge constant of \eqref{eq:Vhatdm}. For continuous weights the
choice of endpoint is immaterial in the limit, as the following corollary shows.

\begin{corollary}[Right-endpoint version]\label{cor:righttag}
Under the conditions of Lemma~\ref{lem:riemann},
\[
\sum_{j=1}^n f_{j}\,(\Delta R_j)^2 \;\plim\; \int_0^T f_u\,d\Lam_u .
\]
\end{corollary}


Under $H_0$, the quadratic covariation $L$ of Lemma~\ref{lem:limitraw}
has an exact conditional distribution given the volatility path, derived
as follows.

\begin{proposition}[Conditional Gaussianity, martingale case]
\label{prop:condgaussraw}
Under Assumptions~\ref{ass:sv}--\ref{ass:mart} and $H_0$,
\[
L\;=\;2\int_0^T h(u)\,dR_u ,\qquad h(u)=\Lam_T-\Lam_u ,
\]
and, conditionally on $\Gsig$,
\[
L\,\big|\,\Gsig \;\sim\; \mathcal{N}\bigl(0,\,V\bigr),
\qquad
V \;=\; 4\int_0^T \bigl(\Lam_T-\Lam_u\bigr)^{2}\,d\Lam_u \;>\;0
\text{ a.s.}
\]
\end{proposition}


The feasible variance estimator replaces $\Lam$ by cumulative realized
variance of the same returns:
\begin{equation}\label{eq:Vhatraw}
\widehat\Lam_n(j)=\sum_{i\le j}r_i^{\,2},
\qquad
\widehat V_n \;=\; 4\sum_{j=1}^{n}
\bigl(\widehat\Lam_n(T)-\widehat\Lam_n(j)\bigr)^{2} r_j^{\,2}.
\end{equation}

\begin{theorem}[Feasible limit theory, martingale case]\label{thm:mart}
Under Assumptions~\ref{ass:sv}--\ref{ass:mart} and $H_0$, as $n\to\infty$,
$\widehat V_n\plim V$ and
\[
t_n\;=\;\frac{\QC_T}{\sqrt{\widehat V_n}}
\;\dlim\;\mathcal{N}(0,1),
\]
the convergence holding conditionally on $\Gsig$ for almost every
volatility path, and hence unconditionally.
\end{theorem}



Theorem~\ref{thm:mart} is a statement about one window of fixed length
$T$. One window information is sufficient for the size of the test under $H_0$. However, under the alternative where $\E[L]\neq0$, it is impossible to estimate $\E[L]$ by one observation window. Hence, to study power we let the number of windows grow and develop an \emph{aggregate mean test}. We partition $[0,KT]$ into the non-overlapping windows
$\mathcal I_k=((k-1)T,\,kT]$, $k=1,\dots,K$, of fixed length $T$, so that
$KT\to\infty$ as $K\to\infty$ while $n\to\infty$ within each window. An
object computed from the returns of window $k$ alone carries the
superscript $(k)$: $\QC^{(k)}$, $\widehat V^{(k)}_n$,
$t^{(k)}_n=\QC^{(k)}/\sqrt{\widehat V^{(k)}_n}$, and the infill limits
$L^{(k)}$, $V^{(k)}$ of Lemma~\ref{lem:limitraw} and
Proposition~\ref{prop:condgaussraw}, with $R^{(k)}_u=X_u-X_{(k-1)T}$ and
$\Lam$ restarted at $(k-1)T$. Set
\[
\tau^{(k)}:=\frac{L^{(k)}}{\sqrt{V^{(k)}}},\qquad
\bar t_{K,n}:=\frac1K\sum_{k=1}^{K}t^{(k)}_n,\qquad
\bar\tau_K:=\frac1K\sum_{k=1}^{K}\tau^{(k)} .
\]
The \emph{aggregate mean test} rejects $H_0$ when
$|\sqrt K\,\bar t_{K,n}|>z_{1-\alpha/2}$.
Assumptions~\ref{ass:sv}--\ref{ass:mart} are imposed on $[0,KT]$. Write $\Gsig_K:=\sigma\bigl((\sigma_u)_{u\le KT}\bigr)$ for the
$\sigma$-field generated by the volatility path over $[0,KT]$, the
$K$-window analogue of $\Gsig$. When imposed, $H_0$ means that
$(\sigma_u)_{u\le KT}$ is independent of $W$. Note that
$t^{(k)}_n\plim\tau^{(k)}$ with or without $H_0$, as proved in Appendix~\ref{app: tstat consistency}.
 
\begin{lemma}[Conditional independence across windows under $H_0$]
\label{lem:condindep}
Under Assumptions~\ref{ass:sv}--\ref{ass:mart} and $H_0$ on $[0,KT]$,
conditionally on $\Gsig_K$ the pairs $(\QC^{(k)},\widehat V^{(k)}_n)$,
$k\le K$, are independent across $k$ for every $n$, and
$(\tau^{(1)},\dots,\tau^{(K)})\mid\Gsig_K\sim\mathcal N(0,I_K)$ a.s.
In particular $\tau^{(1)},\dots,\tau^{(K)}$ are i.i.d.\ $\mathcal N(0,1)$
and independent of $\Gsig_K$.
\end{lemma}
 
 
Under the alternative the $\tau^{(k)}$ are in general neither independent
nor identically distributed, since the conditional location of $L^{(k)}$
is a functional of the volatility path in window $k$ and volatility is
persistent across window boundaries. Power requires only a law of large
numbers, for which the following suffices.
 
\begin{assumption}[Stationary and ergodic window blocks]\label{ass:erg}
The sequence of window blocks
\[
\xi_k:=\Bigl((\sigma_{(k-1)T+u})_{u\in[0,T]},\;
(W_{(k-1)T+u}-W_{(k-1)T})_{u\in[0,T]}\Bigr),\qquad k\ge1,
\]
is strictly stationary and ergodic.
\end{assumption}
 
Since $\tau^{(k)}$ is a fixed measurable functional of $\xi_k$,
$(\tau^{(k)})_{k\ge1}$ is stationary and ergodic. It is integrable under
Assumption~\ref{ass:sv} with or without $H_0$:
$|L^{(k)}|\le2T\overline\sigma^{2}\sup_{u\in\mathcal I_k}|R^{(k)}_u|$ has
all moments by the Burkholder--Davis--Gundy inequality, since
$[R^{(k)}]\le T\overline\sigma^{2}$, and
\begin{equation}\label{eq:Vlower}
V^{(k)}=4\int_{\mathcal I_k}h^{2}\,d\Lam
\;\ge\;4\underline\sigma^{6}\int_0^T(T-u)^{2}du
=\tfrac43\,\underline\sigma^{6}T^{3}\quad\text{a.s.},
\end{equation}
because $h(u)=\Lam_T-\Lam_u\ge\underline\sigma^{2}(T-u)$ and
$d\Lam\ge\underline\sigma^{2}du$. Hence,
\begin{equation}\label{eq:expecttau}
\delta:=\E[\tau^{(1)}]
\end{equation}
is finite.

\begin{corollary}[Size and power of the aggregate mean test,
martingale case]\label{cor:power}
Let Assumptions~\ref{ass:sv}--\ref{ass:mart} hold on $[0,KT]$ and fix
$\alpha\in(0,1)$.
\begin{enumerate}
\item[(i)] \textup{(Size.)} Under $H_0$, for every fixed $K$,
$\sqrt K\,\bar\tau_K\sim\mathcal N(0,1)$ exactly and
$\sqrt K\,\bar t_{K,n}\dlim\mathcal N(0,1)$ as $n\to\infty$. The test has
asymptotic size $\alpha$ for every $K$ with no restriction on the
dynamics of $\sigma$ across windows.
\item[(ii)] \textup{(Power.)} Under Assumption~\ref{ass:erg}, if
$\delta\neq0$ then $\bar\tau_K\to\delta$ a.s.\ and

\[
\lim_{K\to\infty}\liminf_{n\to\infty}
\Prob\!\left(
\left|\sqrt K\,\bar t_{K,n}\right|>z_{1-\alpha/2}
\right)=1.
\]
\end{enumerate}
\end{corollary}
 

\begin{remark}[Scope of the power result]\label{rem:powerscope}
Corollary~\ref{cor:power}(ii) establishes power against $\delta\neq0$,
which is a sub-alternative of $H_1$: by Lemma~\ref{lem:condindep},
$H_0$ implies $\delta=0$, so $\delta\neq0$ implies $H_1$, but the converse
fails. The reason lies in the statistic itself. The quadratic covariation
$[R,R^{2}]$ measures the \emph{directional} asymmetry between returns and
variance innovations: if the dependence between $\sigma$ and $W$ is symmetric in the
sign of the return shock, i.e., $(\sigma,W)\overset{d}{=}(\sigma,-W)$, then
$L$ is odd and $V$ is even in $W$, so $\tau\overset{d}{=}-\tau$ and
$\delta=0$ although $H_0$ fails. 
Such alternatives fall outside the consistency result of Corollary~\ref{cor:power}(ii). A suitably normalized \citet{Neyman1937}'s smooth test based on $J$ Hermite coordinates of the limiting window statistics has an asymptotic $\chi^2_J$ null distribution as $K\to\infty$. This extension can detect additional departures reflected in the selected moments; its feasible finite-sample calibration is left for future research.

Two further considerations put the restriction in perspective. First, for the martingale components of the parametric models used
in Section~\ref{sec:mc}, where the volatility
is an It\^o process and the dependence is carried entirely by the
correlation $\rho$ between the return and variance innovations, one has
$\tau=\rho\,\Upsilon+\sqrt{1-\rho^{2}}\,Z$ with $\Upsilon$ free of $\rho$ and
$Z\sim\mathcal N(0,1)$ independent of $\Upsilon$, so that $\delta(\rho)=a\rho$
with $a=\E[\Upsilon]$ a model constant. Within this family $\delta\neq0$ if and
only if $\rho\neq0$ whenever $a\neq0$, and the alternative of
Corollary~\ref{cor:power} coincides with $H_1$. The regression of
Section~\ref{sec:power} confirms the linearity and estimates $a$.
Second, the economic content of $\delta$ is the leverage effect: according to \eqref{eq:integraCOV}, $\E[L]$ is the integrated covariance between the cumulative return and the spot
variance, of which $\delta$ is the studentized counterpart. The empirical
literature on asymmetric volatility and the leverage effect
\citep{EngleNg1993,AitSahalia2017} documents negative return--variance dependence in equity markets, motivating the directional alternative considered here.
%
The
aggregate mean test should therefore be read as a test of the leverage
effect with independence as its null, rather than as an omnibus test of
dependence.
\end{remark}

\section{Empirical Implementation -- window-constant drift}\label{sec:empimpl}

Assumption~\ref{ass:mart} is restrictive for asset prices: standard asset pricing models predict that assets should have expected returns exceeding the risk free rate which is generally positive. In the presence of a drift, $R_u=A_u+M_u$ with
$A_u=\int_0^u\alpha_s\,ds$ and $M$ the martingale part, the limit of
$\QC_T$ becomes $2\int_0^T(A_u+M_u)\,d\Lam_u$, which carries the
location term $2\int_0^TA_u\,d\Lam_u=O(\|\alpha\|_\infty T\Lam_T)$ --- an
$O(1)$ bias that does \emph{not} vanish as $n\to\infty$. This location
term depends on the unknown drift path and needs not be $\Gsig$-measurable.
Hence, without accounting for it, the martingale-case test $t_n$ of
Theorem~\ref{thm:mart}, which centers at zero, is not valid as it stands.
For the empirical implementation we accordingly assume the drift to be
constant within each window. Demeaning the returns then removes the
location term at its source exactly, path by path
(Lemma~\ref{lem:driftremoval} below), rather than estimating it. This
section develops the demeaned statistic and its inference theory. The proofs of all results in this section, except for the short proof of Lemma~\ref{lem:driftremoval}, are collected
in Appendix~\ref{app:proofs-dm}.

\subsection{The operational statistic: demeaning and exact constant-drift removal}
\label{sec:demean}

We demean the returns within each window,
\[
\rt_j \;=\; r_j-\bar r, \qquad \bar r = \frac1n\sum_{i=1}^n r_i ,
\]
cumulate them into
\begin{equation}\label{eq:bridge-discrete}
\Rt_j \;=\; \sum_{i\le j}\rt_i \;=\; R_j-\frac{j}{n}\,R_T ,
\end{equation}
so that $\Rt_0=\Rt_n=0$, and compute the \emph{demeaned} realized
quadratic covariation
\begin{equation}\label{eq:QCdm}
\QCd_T \;=\; \sum_{j=1}^{n}\Delta\Rt_j\,\Delta\bigl(\Rt^2\bigr)_j
\;=\; \sum_{j=1}^{n}\Bigl\{\,2\,\Rt_{j-1}\bigl(\Delta\Rt_j\bigr)^2
+\bigl(\Delta\Rt_j\bigr)^3\Bigr\}.
\end{equation}
In continuous time the corresponding path is the \emph{bridge} process
$\Rt_u=R_u-(u/T)R_T$, which agrees with \eqref{eq:bridge-discrete} at
every grid point and is pinned to zero at both window endpoints. The point
of the construction is the following exact identity.

\begin{lemma}[Exact removal of within-window-constant drift]\label{lem:driftremoval}
Suppose $\alpha_u\equiv\alpha$ on $[0,T]$ for a (possibly random,
$\F_0$-measurable) constant $\alpha$. Then the demeaned returns computed
from $X$ coincide with the demeaned returns computed from the martingale
part $M$, i.e., $\rt_j=\Delta M_j-\frac1n M_T$ for every $j$, and hence at
every grid point
\[
\Rt_j\;=\;M_j-\frac{t_j}{T}\,M_T .
\]
Consequently every statistic built from demeaned returns, in particular
$\QCd_T$ and $\widehat V^{0}_n$ below, is path by path, exactly equal to
its value computed from the unobserved martingale part.
\end{lemma}

\begin{proof}
$r_j=\alpha\Delta_n+\Delta M_j$ and $\bar r=\alpha\Delta_n+M_T/n$, so
$\rt_j=r_j-\bar r=\Delta M_j-M_T/n$. Then summation gives the result.
\end{proof}

Demeaning therefore reduces the (constant-drift) process to the
martingale case of Section~\ref{sec:mart}, not by making the observed
path a martingale, but by delivering the \emph{bridge} of the martingale.
Having Lemma~\ref{lem:driftremoval}, we are able to do work under
Assumption~\ref{ass:mart} without loss of generality in the remainder of
the section. What remains is to characterize the effect of the bridging
transformation on the limit, its null distribution, and the feasible
variance with the bridge correction. 
%
Bridging reweights the covariance integrand in \eqref{eq:meanQC} by $(1-u/T)$ and preserves the sign of the integrated covariance when the covariance function has a constant sign, as in the Heston benchmark.

For time-varying
drift, demeaning removes the window average $\bar\alpha$, leaving the
centered component $\alpha_u-\bar\alpha$.  In this general case the exact identity of
Lemma~\ref{lem:driftremoval} carries a residual,
\[
\rt_j \;=\; \Delta M_j-\frac1n M_T
\;+\;\int_{t_{j-1}}^{t_j}\bigl(\alpha_u-\bar\alpha\bigr)\,du ,
\qquad
\bar\alpha=\frac1T\int_0^T\alpha_u\,du ,
\]
the added term being the centered-drift increment over the $j$-th
interval, which vanishes when $\alpha$ is constant and recovers
Lemma~\ref{lem:driftremoval}. For the It\^o drift of
a Heston-type stochastic-volatility model its contribution to the expected quadratic covariation vanishes exactly at stationarity
(Appendix~\ref{app:voldrift}).
Its full pathwise effect is also characterized and benchmarked against the window-constant drift assumption (see details in Remark~\ref{rem:hestondrift} and
Table~\ref{tab:driftcorr}).

\subsection{Inference under demeaning: the bridge correction}
\label{sec:bridgeinf}

\begin{lemma}[Infill limit, demeaned statistic]\label{lem:limitdm}
Under Assumption~\ref{ass:sv}, as $n\to\infty$,
\begin{equation}\label{eq:limit}
\QCd_T \;\plim\; L^{0} \;:=\; \bigl[\Rt,\Rt^2\bigr]_T
\;=\; 2\int_0^T \Rt_u\, d\Lam_u .
\end{equation}
\end{lemma}

\begin{remark}[Relation to classical realized-covariation theory]
\label{rem:bnsrelation}
For the undemeaned statistic, the convergence is classical
(Lemma~\ref{lem:limitraw}). Lemma~\ref{lem:limitdm} nonetheless cannot be obtained by direct citation of the classical result. The reason why the extension is needed is that the bridge process
$\Rt_u=R_u-(u/T)R_T$ anticipates the terminal value $R_T$, so it is not
adapted to $(\F_u)$ and is not a semimartingale in that filtration.
The proof of Lemma~\ref{lem:limitdm} (Appendix~\ref{prf:limitdm}) circumvents this by factoring the
anticipative component as the fixed random scalar $R_T$ multiplying
adapted Riemann sums, to which the classical convergence of weighted
realized quadratic variation \citep{jacod1987limit} does apply
(multiplication by the fixed scalar preserving the limit by continuous
mapping).
\end{remark}

\begin{proposition}[Conditional Gaussianity and the bridge-corrected
variance]\label{prop:condgaussdm}
Under Assumption~\ref{ass:sv}, a window-constant drift, and $H_0$,
\begin{equation}\label{eq:rep}
L^{0} \;=\; 2\int_0^T h^{0}(u)\,dR_u ,
\qquad
h^{0}(u) \;=\; \bigl(\Lam_T-\Lam_u\bigr)-c \;=\; h(u)-c,
\qquad
c \;=\; \frac1T\int_0^T u\,d\Lam_u ,
\end{equation}
and, conditionally on $\Gsig$,
\begin{equation}\label{eq:V}
L^{0}\,\big|\,\Gsig \;\sim\; \mathcal{N}\bigl(0,\,V^{0}\bigr),
\qquad
V^{0} \;=\; 4\int_0^T h^{0}(u)^2\,d\Lam_u
\;=\;4\int_0^T\bigl\{(\Lam_T-\Lam_u)-c\bigr\}^2\,d\Lam_u .
\end{equation}
Moreover $V^{0}>0$ a.s.
\end{proposition}


\begin{remark}[The role of the bridge term $c$]\label{rem:bridge}
Relative to the martingale case, bridging replaces the weight
$h(u)=\Lam_T-\Lam_u$ of Proposition~\ref{prop:condgaussraw} by its
centered version $h^{0}=h-c$ and \emph{reduces} the null
variance. Taking a constant volatility $\sigma$ for instance,
\[
V=\tfrac43\sigma^6T^3 ,
\qquad
V^{0}=4\int_0^T\Bigl\{\sigma^2(T-u)-\tfrac{\sigma^2T}{2}\Bigr\}^2\sigma^2\,du
=\tfrac13\sigma^6T^3 ,
\]
so the variance has a fourfold reduction. Omitting $c$ when returns are demeaned therefore yields a
conservative test with standard errors twice large in the
constant-volatility benchmark.
\end{remark}

The feasible estimator replaces $\Lam$ and $c$ by their realized
counterparts computed from the \emph{demeaned} returns:
\begin{equation}\label{eq:Vhatdm}
\widehat\Lam^{0}_n(j)=\sum_{i\le j}\rt_i^{\,2},
\qquad
\widehat c=\sum_{j=1}^n \frac{j}{n}\,\rt_j^{\,2},
\qquad
\widehat V^{0}_n \;=\; 4\sum_{j=1}^{n}
\Bigl\{\bigl(\widehat\Lam^{0}_n(T)-\widehat\Lam^{0}_n(j)\bigr)-\widehat c\Bigr\}^2
\rt_j^{\,2} .
\end{equation}

\begin{theorem}[Consistency of the bridge-corrected variance
estimator]\label{thm:consistency}
Under Assumption~\ref{ass:sv}, a window-constant drift, and $H_0$, as
$n\to\infty$, $\widehat V^{0}_n \plim V^{0}$.
\end{theorem}


\begin{theorem}[Asymptotic standard normality under the null]\label{thm:tstat}
Let Assumption~\ref{ass:sv}, a window-constant drift, and $H_0$ hold, and
define $t^{0}_n = \QCd_T/\sqrt{\widehat V^{0}_n}$. Then, as $n\to\infty$,
\[
t^{0}_n \;\dlim\; \mathcal{N}(0,1),
\]
and the convergence holds conditionally on $\Gsig$ for almost every
volatility path, and hence unconditionally.
\end{theorem}


\begin{remark}[Mixed Gaussianity and pivotality]
Conditionally on $\Gsig$, $L^{0}$ is Gaussian with the deterministic
variance $V^{0}$ of Proposition~\ref{prop:condgaussdm}. But $V^{0}$ is
itself a functional of the volatility path, so the \emph{unconditional}
law of $L^{0}$ is a variance mixture of normals over volatility paths,
with heavier-than-Gaussian tails, the mixed-Gaussian form familiar from \cite{Barndorff-N&Shephard2004a}.
This path dependence is exactly what studentization removes. Because
$V^{0}$ is $\Gsig$-measurable, dividing by $\sqrt{\widehat V^{0}_n}$, which
estimates $V^{0}$ path by path, cancels the path-specific variance before
the mixture over paths is taken: $L^{0}/\sqrt{V^{0}}$ is exactly standard
normal on every path, and the mixture collapses. The studentized statistic
is therefore asymptotically \emph{pivotal}, i.e., its limit law is standard
normal regardless of the volatility dynamics (e.g., level, persistence, long
memory). The nuisance has been divided out rather
than modelled.
\end{remark}

\subsection{Power by aggregation across windows}
Corollary~\ref{cor:power} transfers to the operational statistic, because
demeaning reduces each window to its martingale part. With the notation
preceding Lemma~\ref{lem:condindep}, let $\widehat{QC}^{\,0,(k)}$,
$\widehat V^{0,(k)}_n$,
$t^{0,(k)}_n=\widehat{QC}^{\,0,(k)}/\sqrt{\widehat V^{0,(k)}_n}$,
$L^{0,(k)}$ and $V^{0,(k)}$ be computed from the demeaned returns of
window $k$ as in Section~\ref{sec:demean}, and set
$\tau^{0,(k)}:=L^{0,(k)}/\sqrt{V^{0,(k)}}$,
$\bar t^{\,0}_{K,n}:=K^{-1}\sum_kt^{0,(k)}_n$,
$\bar\tau^{\,0}_K:=K^{-1}\sum_k\tau^{0,(k)}$. The drift is window-constant,
$\alpha_u\equiv\alpha^{(k)}$ on $\mathcal I_k$ with $\alpha^{(k)}$
$\F_{(k-1)T}$-measurable. The proof of Lemma~\ref{lem:driftremoval} applies
\emph{mutatis mutandis} to each window.
 
\begin{corollary}[Aggregate mean test, demeaned statistic]\label{cor:powerdm}
Let Assumption~\ref{ass:sv} and a window-constant drift hold on $[0,KT]$.
Then Lemma~\ref{lem:condindep} and Corollary~\ref{cor:power} hold with
\begin{align*}
&\bigl(\widehat{QC}^{\,0,(k)},\widehat V^{0,(k)}_n,t^{0,(k)}_n,\tau^{0,(k)},
  \bar t^{\,0}_{K,n},\bar\tau^{\,0}_K\bigr)
&\textup{in place of}\quad
\bigl(\QC^{(k)},\widehat V^{(k)}_n,t^{(k)}_n,\tau^{(k)},\bar t_{K,n},\bar\tau_K\bigr),
\end{align*}
Assumption~\ref{ass:mart} dropped, and \eqref{eq:Vlower} replaced by
$V^{0,(k)}\ge\tfrac13\underline\sigma^{6}T^{3}$ a.s.
\end{corollary}
 
\noindent The proof of Corollary~\ref{cor:powerdm} is given in Appendix~\ref{prf:powerdm}.

\begin{remark}[On the window-constant-drift assumption]\label{rem:constdrift}
The window-constant drift is an assumption of convenience: it is what
makes Lemma~\ref{lem:driftremoval} exact and lets every result of this
section reduce to the martingale case of Section~\ref{sec:mart}. It is
not, however, a property one can verify on the data, and standard
continuous-time models violate it. In any stochastic-volatility model
with constant proportional price drift the log-return carries the It\^o
term $-\tfrac12v_u$ (from Equation \eqref{XDGP}), which fluctuates with the variance path within
every window, so the true drift is window-constant only if the volatility
is. Two considerations justify working under the assumption regardless.
First, Section~\ref{sec:demean} showed that demeaning removes the window
average $\bar\alpha$ of any drift, leaving only the centered residual
$\alpha_u-\bar\alpha$. Over a monthly window this residual is small
relative to the diffusive increments, and the smaller the window the
smaller it is. Second, for the Heston model of Section~\ref{sec:mc}, the
effect of the residual can be computed exactly rather than bounded.
Appendix~\ref{app:voldrift} shows that the infill limit of the demeaned
statistic decomposes pathwise as $L^{0}=L^{0}_{M}+D$, where $L^{0}_{M}$
is the martingale-part limit to which Proposition~\ref{prop:condgaussdm}
applies and $D$ is a $\Gsig$-measurable functional of the variance path
alone. Conditionally on the volatility one can compute a feasible correction $\widehat D$ that restores exact asymptotic normality
without the assumption. 
Further discussion and impact quantification are shown in Remark~\ref{rem:hestondrift} and
Table~\ref{tab:driftcorr}.
\end{remark}

\begin{remark}[Jumps]\label{rem:jumps}
The \cite{LeeMyklandJump} detector flags a return as a jump when its
studentized magnitude exceeds a threshold $\beta_n$, the percentile of
the limiting law of the maximum of the detection statistic under the
diffusive null. As $\Delta_n\to0$ with $\beta_n\to\infty$ slower than
$\sqrt{n\log n}$, both global error rates vanish, i.e., the probability of
failing to detect an actual jump and the probability of spurious
detection, the latter at rate $\exp(-\beta_n)$
\citep[Theorems 3--5]{LeeMyklandJump}. All finite-activity jumps are
therefore classified correctly with probability tending to one, so
conditioning on correct classification is asymptotically justified. This asymptotic result concerns a diverging threshold. The simulations use the fixed level $\alpha_{\mathrm{LM}}=0.01$ and assess the finite-sample performance of that implementation. At any
fixed frequency the two error types persist and respond oppositely to the
threshold, raising $\beta_n$ suppresses spurious zeroing, which biases
$\widehat\Lam^{0}$ (hence $\widehat V^{0}_n$) downward, but admits more
undetected jumps, which distort the distribution of $\QCd_T$ (a jump
contributes $\approx 2\Rt_{\tau-}J^2+J^3$). This residual contamination is
not assumed away: it is present in the Monte Carlo of
Section~\ref{sec:mc}, which applies the detector to simulated jump paths,
and is reflected in the reported sizes.
\end{remark}

\section{Simulation Study}\label{sec:mc}

In order to heuristically validate this theory and to measure the power of the developed tests, we perform a simulation experiment, which generates multiple paths of a window of high-frequency observations, each consists of $22$ trading days with $78$ five-minute returns per day. These simulation paths are derived by the following data-generating process (DGP), which specializes our general semimartingale \eqref{eq:smdecomp}:
\begin{align}
	\label{SDGP}
	\mathrm{d}S_t &= \left(
	\mu - \bar k\right)S_t\, \mathrm{d}t + \sqrt{v_t} S_t\,\mathrm{d}W^x_t + S_tJ_t\mathrm{d}N_t , \\
		\label{VDGP}
	\mathrm{d}v_t &= \kappa(\theta - v_t)\mathrm{d}t + \xi \sqrt{v_t} \mathrm{d}W^v_t, \\
	\label{corr}
	&\mathbb{E}[\mathrm{d}W^x_t \, \mathrm{d}W^v_t] = \rho \mathrm{d}t.
\end{align}
Equation~\eqref{SDGP} specifies the equity price $S_t$ with a continuous diffusion component of drift $\mu$ and stochastic variance $v_t$, augmented by multiplicative jumps of relative size $J_t$ arriving according to a homogeneous Poisson process $N_t$ with intensity $\lambda \in \mathbb{R}_+$. 
We specify multiplicative jump sizes by $\log(1+J_t)\sim \mathcal{N}\!\left(\log(1+\mu_J)-\sigma_J^2,\,2\sigma_J^2\right)$, so that $\E[J_t]=\mu_J$.
Following \citet{bates1996jumps}, the compensator $\bar{k}=\lambda \mu_J$ offsets the average jump effect so that $\mathbb{E}[\mathrm{d}S_t \mid S_t]=\mu S_t\,\mathrm{d}t$. In~\eqref{VDGP}, the instantaneous variance follows the CIR process \citep{cox1985theory}, where $\kappa$ is the speed of mean reversion, $\theta$ the long-run mean, and $\xi$ the volatility of variance. The Brownian motions $W_t^x$ and $W_t^v$ are correlated as in~\eqref{corr}, with $\rho \in [-1,1]$.
We treat $(\mu,\kappa,\theta,\xi)$ as constants. For the CIR dynamics, a strong solution exists for any $(\kappa,\theta,\xi)\in(0,\infty)^3$ and $v_0>0$ \citep{RevuzYor1999}. If $\gamma \coloneqq 2\kappa\theta/\xi^2 \ge 1$ (the Feller condition), $v_t$ remains strictly positive; when $0\le\gamma<1$, the origin is attainable and instantaneously reflecting \citep{Feller1951,PitmanYor1982}. We assume $\kappa,\theta,\xi>0$, and the Feller condition satisfied.

For affine tractability \citep{duffie2000transform}, define $X_t=\log S_t$. By It\^o's lemma \citep{ORS2022}, we obtain the dynamic:
\begin{align}
\label{XDGP}
\mathrm{d}X_t \;=\; \Bigl(\mu - \tfrac{1}{2} v_t - \tilde{k}\Bigr)\,\mathrm{d}t + \sqrt{v_t}\,\mathrm{d}W_t^x + \log(1+J_t)\,\mathrm{d}N_t,
\end{align}
with $\tilde{k} \coloneqq \lambda\!\left[\log(1+\mu_J)-\sigma_J^2\right]$ compensating the average jump contribution in log space. Equations \eqref{VDGP}, \eqref{corr}, and \eqref{XDGP} now jointly define the system of stochastic differential equations (SDEs) governing our data-generating process. We call this system the Bates' framework. In this framework, the parameter $\rho$, which parametrizes the dependence between the return and variance innovations, i.e., $W^x$ and $W^v$, plays an important role since it generates the quadratic covariation the test targets, i.e., $\rho=0$ is the $H_0$. With $\lambda=0$ the system reduces to the Heston model \citep{heston1993closed}. 

Diffusive paths are generated by the quadratic-exponential (QE) discretization of
\citet{andersen2008simple} while jumps are overlaid
by drawing $N_T\sim\mathrm{Poisson}(\lambda T)$ with arrival times uniform
on $[0,T]$ and sizes drawn as above. We set $v_0=\theta$, initializing variance at its long-run mean.
The baseline setting is
$\kappa=1.5$, $\theta=0.04$, $\xi=0.3$ (Feller ratio $2\kappa\theta/\xi^2=1.33$), and $\mu=\theta/2$, so that
the expected diffusive drift $\mu-\E[v_t]/2$ is zero. The jump-free
experiments ($\lambda=0$) run six configurations, $32{,}000$ replications
each for A--F as shown in Table \ref{tab:jumpfree}.
\begin{table}[htbp]
\captionsetup{justification=raggedright,singlelinecheck=false}
\caption{The jump-free experiment settings.}
\centering
\label{tab:jumpfree}
\begin{tabular}{lll}
\toprule
Config & Other parameters & Verification target \\
\midrule
A\; $\rho=0$    & Baseline & Size \\
B\; $\rho=-0.3$ & Baseline & Power \\
C\; $\rho=-0.7$ & Baseline & Power \\
D\; $\rho=-0.9$ & Baseline & Power \\
E\; $\rho=0$    & Baseline except $\mu=\theta/2+0.10$
  & Size with demeaning robustness \\
F\; $\rho=0$    & $\theta\sim U[0.05,0.5]$, $\xi\sim U[0.1,1]$,
  $\gamma\sim 1+\mathrm{Exp}(4)$, & Size uniformity over \\
  &  $\kappa=\gamma\xi^2/(2\theta)$, $\mu=\theta/2$
  & parameter uncertainty \\
\bottomrule
\end{tabular}
\end{table}
The jump experiments of Section~\ref{sec:batessize} vary
$(\lambda,\mu_J,\sigma_J)$ around this baseline. For each simulated window
we compute $\QCd_T$ from demeaned returns and studentize it in two ways: (i)
by using the BNS-type variance estimator $H_n$ of \eqref{eq:Hn}; (ii) by using the
conditional-Gaussian estimator $\widehat V^{0}_n$ of \eqref{eq:Vhatdm},
including the bridge term $\widehat c$. We then run a two-sided hypothesis test on the studentized statistics for each window at the 5\% significance level.

\subsection{Heston Framework (jump-free)}
\subsubsection{Quadratic Covariation in a closed form}\label{sec:hestonbench}

The expected quadratic covariation in Heston framework admits a closed-form expression, allowing the simulator and the theory to validate one another. The derivation is given in Appendix~\ref{app:heston} and the following summarizes our results.

\begin{proposition}[Mean of the quadratic covariation under Heston
dynamics with a window-constant drift]\label{prop:heston}
Let the variance follow \eqref{VDGP} with $\lambda=0$ and $v_0=\theta$,
and let the log-return follow \eqref{XDGP} with its drift
$\mu-\tfrac12v_u$ replaced by a window-constant $\alpha$. Let $\phi=\kappa T$. 
Then the return--variance covariance function is
\begin{equation}\label{eq:mheston}
m(u)\;:=\;\E\bigl[R_u\,v_u\bigr]
\;=\;\frac{\rho\,\xi\,\theta}{\kappa}\bigl(1-e^{-\kappa u}\bigr).
\end{equation}
By Lemma~\ref{lem:driftremoval} the demeaned statistic is pathwise identical
to the one computed from the martingale part $dR_u=\sqrt{v_u}\,dW^x_u$, on
which the expectations below are evaluated. The expected quadratic covariation without demeaning refers to that martingale part directly. 
Therefore, the infill limits $L$ and $L^{0}$ of the undemeaned and demeaned
statistics (Lemmas~\ref{lem:limitraw} and~\ref{lem:limitdm}) have expectations
\begin{align}
\E\bigl[L\bigr]
&=\frac{2\rho\xi\theta}{\kappa}
\left(T-\frac{1-e^{-\phi}}{\kappa}\right),
\label{eq:qcovraw}\\[2pt]
\E\bigl[L^{0}\bigr]
&=\frac{2\rho\xi\theta}{\kappa}
\left(\frac{T}{2}-\frac{1}{\kappa}+\frac{1-e^{-\phi}}{\kappa^{2}T}\right).
\label{eq:qcovdm}
\end{align}
Both brackets are strictly positive for every $\kappa T>0$. Consequently
both means vanish if and only if $\rho=0$.
\end{proposition}

Replacing the It\^o drift $-\tfrac12v_u$ by a constant is what makes
Lemma~\ref{lem:driftremoval} applicable and the proof elementary. The
simulations of Section~\ref{sec:size} use the exact dynamics
\eqref{XDGP}. Appendix~\ref{app:voldrift} shows that the omitted term
enters $L^{0}$ only through a $\rho$-independent, $\Gsig$-measurable
shift $D$, whose effect is quantified in Table~\ref{tab:driftcorr}.

Four features of \eqref{eq:qcovraw}--\eqref{eq:qcovdm} organize the
experiments. (i) \emph{Identification:} the expected quadratic covariation is exactly linear in $\rho$ and proportional to $\xi\theta$, vanishing precisely at
$\rho=0$, i.e., the closed-form version of the leverage interpretation of
Section~\ref{sec:svsm} and the equivalence discussed around \eqref{eq:meanQC}.
(ii) \emph{Persistence weighting:} in the slow-mean-reversion case, 
as $\kappa T\to0$, the $\E[L]$ and $\E[L^{0}]$ converge to $\rho\xi\theta T^{2}$ (raw) and $\rho\xi\theta T^{2}/3$ (demeaned), while for $\kappa T\to\infty$
the limits are $2\rho\xi\theta T/\kappa$ and $\rho\xi\theta T/\kappa$.
The quadratic covariation aggregates leverage weighted by the durability of the
volatility response, not the instantaneous $\rho$ alone. (iii) \emph{Bridge
cost:} demeaning shrinks the $\E[L]$ to $\E[L^{0}]$ by a factor between $3$ (slow
reversion) and $2$ (fast reversion), just as it shrinks the null variance to $V^{0}$ in Remark~\ref{rem:bridge}. So bridging attenuates both the location and the dispersion.
(iv) \emph{Third
cumulant:} the martingale-part third cumulant is
$\kappa_3(R_T)=\tfrac{3}{2}\,\E[L]$
(Appendix~\ref{app:3rdcumulant}), connecting the quadratic covariation to the \cite{NP2020} aggregation, which has been discussed by \cite{choelee}. 

\begin{remark}[Heston as a robustness check on the window-constant-drift assumption]
\label{rem:hestondrift}
The Heston log-price \eqref{XDGP} carries the drift $\mu-\tfrac12 v_u$.
Demeaning removes $\mu$ exactly (Lemma~\ref{lem:driftremoval}), but the
It\^o term $-\tfrac12 v_u$ is not window-constant, so
Proposition~\ref{prop:condgaussdm} and Theorem~\ref{thm:tstat} do not
apply to \eqref{XDGP} literally. The Heston model is, however, tractable
enough to resolve this exactly rather than approximately. In
Appendix~\ref{app:voldrift} we show that the infill limit of the
demeaned statistic decomposes pathwise as
\[
L^{0}\;=\;L^{0}_{M}\;+\;D,
\qquad
D\;=\;\Lam_T\Bigl(c-\tfrac12\Lam_T\Bigr),
\]
where $L^{0}_{M}$ is the demeaned limit of the martingale part, to which
Proposition~\ref{prop:condgaussdm} applies \emph{mutatis mutandis}, and $D$ is a
$\Gsig$-measurable functional of the variance path only. Under $H_0$, conditionally on
$\Gsig$ the stochastic-volatility drift is therefore a pure location
shift: $L^{0}\mid\Gsig\sim\mathcal N(D,V^{0})$, with the same variance
$V^{0}$ as in the martingale case. Two consequences follow. First,
$\E[D]=0$ at stationarity, so the closed form \eqref{eq:qcovdm} is exact.
Second, $D$ has a natural feasible counterpart $\widehat D$ built from the
same demeaned squared returns as $\widehat V^{0}_n$, and the corrected
statistic $t^{0,D}_n=(\QCd_T-\widehat D)/\sqrt{\widehat V^{0}_n}$ is
asymptotically $\mathcal N(0,1)$ under $H_0$ \emph{without} the
window-constant-drift assumption. The uncorrected $t^{0}_n$ of
Theorem~\ref{thm:tstat} differs from it by $\widehat D/\sqrt{\widehat
V^{0}_n}$, which is of order $10^{-2}$ at the calibration of
Table~\ref{tab:driftcorr}. Section~\ref{sec:size} quantifies this on the
simulated paths. Since the true return dynamics are unknown in practice,
$t^{0}_n$ is the statistic we take to the data. The Heston case shows
what that approximation costs when the drift is in fact volatility-driven.
\end{remark}

\subsubsection{Size and validation} \label{sec:size}

\begin{table}[htbp]
\centering
\caption{Jump-free simulation results. Rejections (our method vs. BNS method \citep{Barndorff-N&Shephard2004a}) at two-sided $5\%$ significance.
$t$ denotes the conditional-Gaussian test statistic $t^{0}_n$, whereas $t_{\mathrm{BNS}}$ denotes the BNS test statistic from \eqref{eq:bnst}.
The last two columns compare the Monte Carlo mean of the demeaned realized quadratic covariation with
its infill limit in analytical form from Proposition~\ref{prop:heston}.}
\label{tab:size}
\small
\begin{tabular}{lcccccccc}
\toprule
Config & mean$(t)$ & sd$(t)$ & sd$(t_{\mathrm{BNS}})$ & rej.\ $t$ & rej.\ BNS
& $\E[\QCd]$ MC & $\E[L^{0}]$ analytic \\
\midrule
A\; $\rho=0$                  & $+0.003$ & $0.999$ & $17.70$ & $0.049$        & $0.955$        & $+3.49\text{e-}7$ & $0$ \\
B\; $\rho=-0.3$               & $-0.025$ & $0.997$ & $17.65$ & $0.051^{\ast}$ & $0.053^{\ast}$ & $-8.46\text{e-}6$ & $-8.08\text{e-}6$ \\
C\; $\rho=-0.7$               & $-0.052$ & $1.000$ & $17.66$ & $0.053^{\ast}$ & $0.059^{\ast}$ & $-1.88\text{e-}5$ & $-1.89\text{e-}5$ \\
D\; $\rho=-0.9$               & $-0.065$ & $0.994$ & $17.68$ & $0.051^{\ast}$ & $0.062^{\ast}$ & $-2.42\text{e-}5$ & $-2.42\text{e-}5$ \\
E\; $\rho=0$, drift $+10\%$/y & $+0.001$ & $0.997$ & $17.66$ & $0.050$        & $0.953$        & $-2.45\text{e-}7$ & $0$ \\
F\; random parameters         & $-0.003$ & $0.998$ & $17.69$ & $0.049$        & $0.954$        & --- & --- \\
\bottomrule
\end{tabular}

\vspace{2pt}
{\footnotesize\raggedright
$^{\ast}$ Size-adjusted rejection rates. For the power configurations
(B--D) the entries in columns rej.\ $t$ and rej.\ BNS use size-adjusted
critical values instead of $\pm1.96$: the empirical $2.5$th and $97.5$th
percentiles of each $t$-statistic under the size configuration~A ($\rho=0$)
are taken as critical values, and the rejection rate is the fraction
of the $t$-statistic falling outside them. This isolates power net of any
size distortion, which matters for BNS: its raw rejection rate is
$\approx0.95$ irrespective of $\rho$ ($t_{\mathrm{BNS}}$ diverges, with
$\mathrm{sd}(t_{\mathrm{BNS}})\approx17.7$), so the raw rejection rate is
uninformative about asymmetry, whereas the size-adjusted rejection reveals
whether either $t$-statistic detects $\rho\neq0$ per window. Configurations
A, E, and F have $\rho = 0$ and retain the $5\%$ $\pm1.96$ critical
values.\par}
\end{table}

Table~\ref{tab:size} reports the jump-free simulation results. 
The conditional-Gaussian statistic attains the
correct size in the null configurations (A, E, and F): rejection rates
between $4.9\%$ and $5.0\%$ at the $5\%$ significance level, with mean close to zero and
standard deviation close to one. 
These observations demonstrate accurate finite-sample calibration in the Heston experiments, including the volatility-dependent drift examined in Table~\ref{tab:driftcorr}.
%
%
%
The BNS statistic centered at zero rejects in approximately 95\% of the independence-null simulations, consistent with the distinction between hypotheses discussed in Section~\ref{sec:twonulls}.

Second, under the alternatives (B--D) neither $t$-statistic has useful
per-window power: the size-adjusted rejection rates are low for both, because the per-window location shift is small relative to the
dispersion of $L^{0}$.
This is not a deficiency of either test but a property of the estimand, i.e., at a one-month horizon a single window carries little information about
the asymmetry. Power is therefore obtained by aggregating across windows rather than
within a single window. By Corollary~\ref{cor:powerdm}, under $H_0$, the
aggregate mean test based on $\sqrt K\,\bar t^{\,0}_{K,n}$ has asymptotic
size $\alpha$ for every $K$ while under the alternative, power tends to
one as $K\to\infty$ with the window length $T$ fixed (under Assumption~\ref{ass:erg}). Hence, with the window length fixed, power can be recovered by aggregating the standardized statistics across windows, which we investigate empirically next. 
Third, the Monte Carlo means of $\QCd_T$ reproduce the closed form
\eqref{eq:qcovdm} across all leverage strengths, validating simulator and
theory jointly.

\begin{table}[H]
\centering
\caption{Impact of the stochastic-volatility drift correction (Appendix~\ref{app:voldrift}).
For each configuration of Table~\ref{tab:size}, $t^{0}_n$ is the conditional-Gaussian statistic of
Theorem~\ref{thm:tstat} and $t^{0,D}_n$ is its drift-corrected counterpart \eqref{eq:Dhat}.
Columns 2--4 summarize the per-window difference $t^{0}_n-t^{0,D}_n=\widehat D/\sqrt{\widehat V^{0}_n}$
across replications. $\E[\widehat D]$ MC is the Monte Carlo mean of the feasible drift term.
$\mathrm{corr}(t^{0}_n,t^{0,D}_n)$ is the correlation between
the two statistics. The last two columns are their two-sided rejection rates at 5\% significance.}
\label{tab:driftcorr}
\small
\begin{tabular}{lcccccccc}
\toprule
& \multicolumn{3}{c}{$t^{0}_n-t^{0,D}_n$} & & & \multicolumn{2}{c}{rej.\ at $\pm1.96$} \\
\cmidrule(lr){2-4}\cmidrule(lr){7-8}
Config & mean & sd & $\max|\cdot|$ & $\E[\widehat D]$ MC & $\mathrm{corr}(t^{0}_n,t^{0,D}_n)$
& $t^{0}_n$ & $t^{0,D}_n$ \\
\midrule
A\; $\rho=0$                  & $-2.9\text{e-}4$ & $3.8\text{e-}3$ & $1.3\text{e-}2$ & $+7.8\text{e-}8$ & $0.999993$ & $0.049$ & $0.049$ \\
B\; $\rho=-0.3$               & $-3.1\text{e-}4$ & $3.7\text{e-}3$ & $1.3\text{e-}2$ & $+7.5\text{e-}8$ & $0.999993$ & $0.049$ & $0.049$ \\
C\; $\rho=-0.7$               & $-3.1\text{e-}4$ & $3.8\text{e-}3$ & $1.7\text{e-}2$ & $+7.8\text{e-}8$ & $0.999993$ & $0.051$ & $0.051$ \\
D\; $\rho=-0.9$               & $-2.7\text{e-}4$ & $3.7\text{e-}3$ & $1.5\text{e-}2$ & $+7.9\text{e-}8$ & $0.999993$ & $0.049$ & $0.049$ \\
E\; $\rho=0$, drift $+10\%$/y & $-3.3\text{e-}4$ & $3.8\text{e-}3$ & $1.8\text{e-}2$ & $+7.4\text{e-}8$ & $0.999993$ & $0.050$ & $0.050$ \\
F\; random parameters         & $-4.5\text{e-}4$ & $7.7\text{e-}3$ & $4.1\text{e-}2$ & $+2.4\text{e-}6$ & $0.999971$ & $0.049$ & $0.049$ \\
\bottomrule
\end{tabular}
\end{table}

Table~\ref{tab:driftcorr} isolates the effect of the stochastic-volatility
drift correction of Appendix~\ref{app:voldrift}. Across all configurations
the per-window difference $t^{0}_n-t^{0,D}_n=\widehat D/\sqrt{\widehat V^{0}_n}$
has a Monte Carlo mean of order $10^{-4}$ and never exceeds $2\times10^{-2}$
on the $\mathcal N(0,1)$ scale (and $4\times10^{-2}$ under the
random-parameter draws of config~F), while the Monte Carlo mean of $\widehat D$ is small under the deterministic initialization $v_0=\theta$. The uncorrected and corrected statistics
are therefore virtually collinear, with $\mathrm{corr}(t^{0}_n,t^{0,D}_n)>0.9999$,
and their rejection rates coincide to within Monte Carlo error. The window-constant-drift approximation has a small effect in the reported monthly calibrations.

\subsubsection{Power and its interpretation} \label{sec:power}
Under the alternative, by Corollary~\ref{cor:powerdm}, the
aggregate mean test is consistent whenever the per-window shift
$\delta:=\E[\tau^{0,(1)}]$ is nonzero: $\bar\tau^{\,0}_K\to\delta$ a.s., so
$\sqrt K\,\bar\tau^{\,0}_K\to\pm\infty$ and
$\Prob\bigl(|\sqrt K\,\bar t^{\,0}_{K,n}|>z_{0.975}\bigr)\to1$ as first
$n\to\infty$ and then $K\to\infty$. The present section quantifies it. 
As we have seen, a single window is individually too small to detect asymmetry (Table~\ref{tab:size}). Aggregating $K$ windows amplifies the location by $\sqrt K$ against an $O(1)$ dispersion and thereby recovers power. In finite samples the per-window shift $\delta$ (Equation~\eqref{eq:expecttau}) is estimated by the sample average $\hat\delta$ in Monte Carlo, which enters
the analytic power below.

Then, based on the jump-free experiment, in particular configurations B--D, we
apply a resampling technique to recover the power, in the following steps. 

\paragraph{Setup.} Fix the value of leverage parameter $\rho$ in the selected configurations of the jump-free experiment (configurations B--D under the alternative,
configuration~A under the null), and consider the pool of $B=32{,}000$
per-window statistics $\{t^{0}_n\}$ simulated under it, one value per
window. The independently simulated windows provide an i.i.d. benchmark
that isolates the effect of aggregation on power.
 
\paragraph{Step 1 (aggregation).} Draw $K$ statistics with replacement from
the pool and compute $\sqrt K\,\bar t^{\,0}_{K,n}$, where $\bar t^{\,0}_{K,n}$
is the sample average. Each $\bar t^{\,0}_{K,n}$ is also recentered by the following step.
 
\paragraph{Step 2 (recentering).} In a finite pool, the mean of $t^{0}_n$
under $H_0$ is not exactly zero: configuration~A gives
mean$(t)=+0.003$ (Table~\ref{tab:size}). Because resampling treats
the pool as the population, this deviation is amplified by $\sqrt K$ exactly
like a genuine shift, e.g.\ to $0.003\sqrt{24{,}000}\approx0.46$ at the
largest $K$, and would distort the rejection rate even in the absence of
any true leverage effect. We therefore replace each $t^{0}_n$ by $t^{0}_n-\bar t_{\mathrm A}$,
where $\bar t_{\mathrm A}$ is the configuration~A pool mean, so that power is evaluated relative to the simulated null reference. 
 
\paragraph{Step 3 (simulated power).} Reject $H_0$ when the recentered
aggregate satisfies $|\sqrt K\,\bar t^{\,0}_{K,n}|>z_{0.975}$. Repeating
Steps~1--2 over $16{,}000$ bootstrap replications and recording the fraction
of rejections gives the resampling-based power at that $\rho$ and $K$,
reported in Table~\ref{tab:aggpower}.
 
\paragraph{Step 4 (analytic power).} We also consider an analytic power by the Gaussian approximation for benchmark comparison. If the Gaussian approximation
$\sqrt K\,(\bar t^{\,0}_{K,n}-\delta)\approx\mathcal N(0,1)$ is accurate, the
two-sided rejection probability is available in closed form,
\[
\mathrm{Power}(\rho;K)=\Phi\bigl(-z_{0.975}+|\delta|\sqrt{K}\bigr)
+\Phi\bigl(-z_{0.975}-|\delta|\sqrt{K}\bigr),
\]
without any resampling. In practice $\delta$ is unknown and is replaced by
$\hat\delta$, the mean of the $B$ simulated statistics at that $\rho$,
recentered as in Step~2.

\begin{table}[htbp]
\centering
\caption{Power of the aggregate mean test $\sqrt K\,\bar t^{\,0}_{K,n}$ (conditional-Gaussian
statistic) at the two-sided $5\%$ significance level. Each entry is the bootstrap
rejection frequency over $16{,}000$ draws from a pool of $32{,}000$ simulated
windows. Config~A ($\rho=0$) is the size row. $K$ denotes the number of resamples
used to compute $\sqrt K\,\bar t^{\,0}_{K,n}$.}
\label{tab:aggpower}
\small
\begin{tabular}{lcccccc}
\toprule
Config & $K=60$ & $K=120$ & $K=240$ & $K=1{,}200$ & $K=6{,}000$ & $K=24{,}000$ \\
\midrule
A\; $\rho=0$    & $0.050$ & $0.050$ & $0.049$ & $0.048$ & $0.049$ & $0.051$ \\
B\; $\rho=-0.3$ & $0.055$ & $0.058$ & $0.072$ & $0.164$ & $0.591$ & $0.992$ \\
C\; $\rho=-0.7$ & $0.072$ & $0.092$ & $0.138$ & $0.477$ & $0.989$ & $1.000$ \\
D\; $\rho=-0.9$ & $0.080$ & $0.117$ & $0.188$ & $0.651$ & $1.000$ & $1.000$ \\
\bottomrule
\end{tabular}
\end{table}

Table~\ref{tab:aggpower} shows near-nominal null rejection rates across all values of $K$ in the recentered independent-window experiment.

Secondly under the three selected alternative cases, the rejection probability converges to one if the resampling number is sufficiently large. 

To have a complete overview of the power across the leverage parameter
space, we trace out the full power curve over $\rho\in(-1,1)$. 
The linear dependence of the expected statistics on $\rho$ motivates the approximation $\delta(\rho)\approx a\rho$ for the simulated statistics. We estimate $a$ by ordinary least squares through the origin.

The regression yields a
single slope $\hat a$, which is estimated from an independent pool of $128{,}000$ simulated windows at each $\rho\in\{-0.9,-0.8,\ldots,0.8,0.9\}$. The fitted analytic power curve at aggregate size
$K$ is then
\[
\mathrm{Power}(\rho;\,K)
\;=\;\Phi\bigl(-z_{0.975}+|\hat a\,\rho|\sqrt{K}\bigr)
\;+\;\Phi\bigl(-z_{0.975}-|\hat a\,\rho|\sqrt{K}\bigr),
\]
i.e.\ the per-configuration $\hat\delta$ of the previous paragraph replaced
by its fitted value $\hat a\,\rho$, and symmetric in $\rho$ by construction.

\begin{figure}[htbp]
\centering
\includegraphics[width=\textwidth]{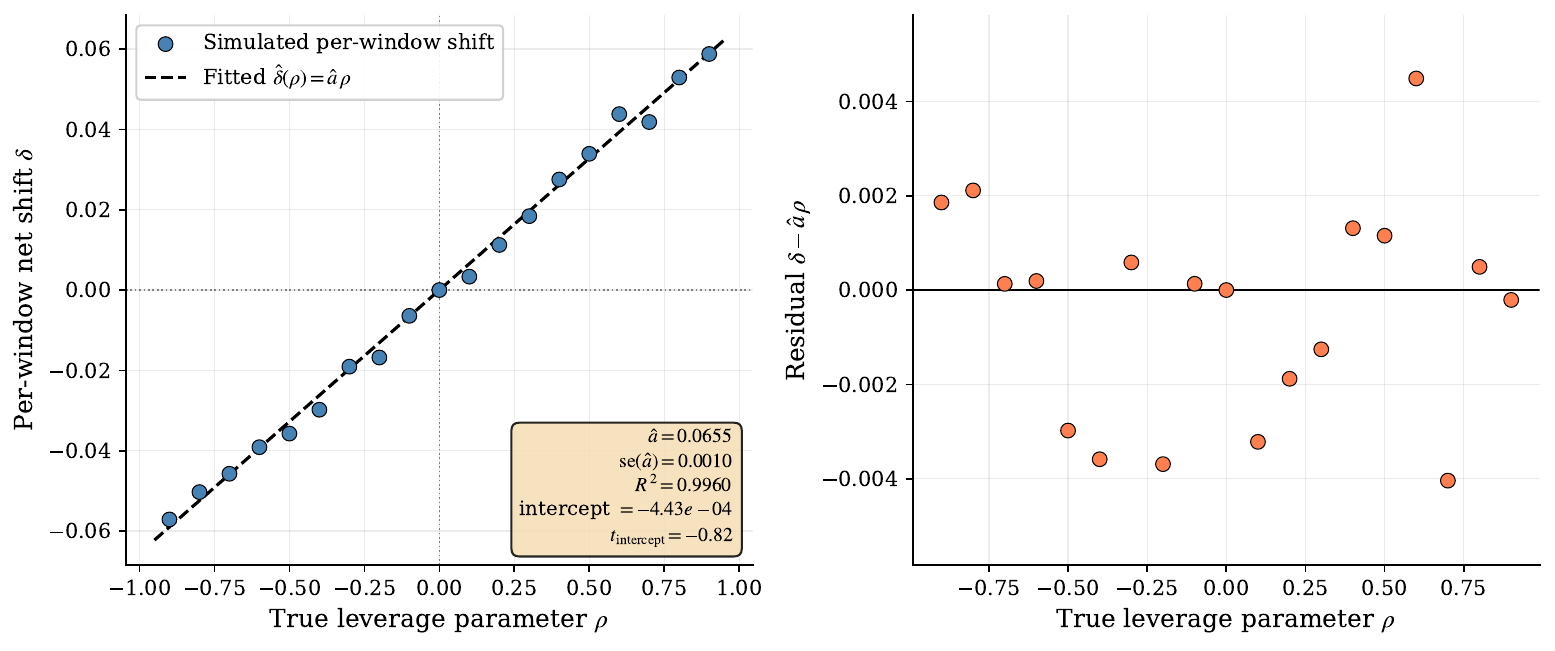}
\captionsetup{justification=raggedright,singlelinecheck=false}
\caption{Regression diagnostic for the per-window shift
$\delta(\rho)=a\,\rho$.}
\label{fig:regression_diagnostic}
\end{figure}

Figure~\ref{fig:regression_diagnostic} confirms the quality of the linear
fit. The left panel plots the simulated per-window shifts (dots) against
the fitted line $\hat\delta(\rho)=\hat a\,\rho$. The points lie tightly close to
the line, with $R^2=0.996$ and a fitted slope
$\hat a=0.0655\;(\mathrm{se}=0.001)$. The unconstrained intercept is
$-4.4\times10^{-4}$ with $|t|=0.82$, indistinguishable from zero and
consistent with the odd-function symmetry
$\delta(-\rho)=-\delta(\rho)$. The right panel displays the residuals,
which scatter randomly around zero with no visible curvature or
heteroskedasticity, confirming that the linear-in-$\rho$ specification
captures the per-window shift across the full leverage range.

\begin{figure}[htbp]
\centering
\includegraphics[width=\textwidth]{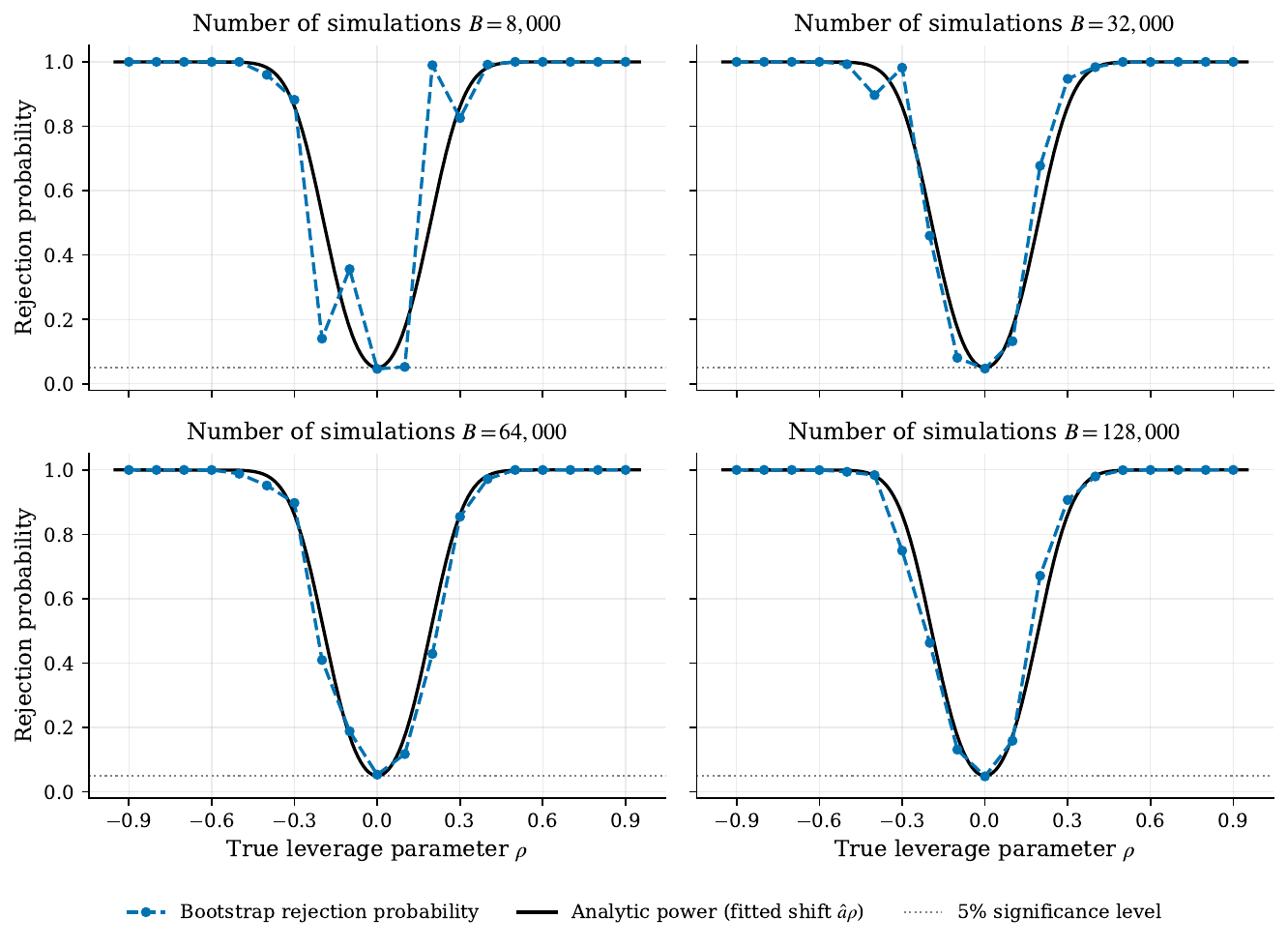}
\caption{Power of the aggregate mean test $\sqrt K\,\bar t^{\,0}_{K,n}$ against the
true leverage parameter $\rho$, with $K=24{,}000$ fixed.}
\label{fig:power_curves}
\end{figure}

Figure~\ref{fig:power_curves} plots the bootstrap rejection probability
(solid) and the regression-fitted analytic power (dashed) against the true
leverage parameter $\rho$, with the aggregate size fixed at $K=24{,}000$.
Four panels correspond to increasing number of simulation paths 
$B\in\{8{,}000,\;32{,}000,\;64{,}000,\;128{,}000\}$, each providing the
pool from which the $K$ resamples are drawn. Three features stand out.
First, the bootstrap power curve converges to the smooth analytic curve as
the pool size $B$ grows: at $B=8{,}000$ the curve is jagged and departs
visibly from the fitted envelope, while at $B=128{,}000$ the two are
nearly indistinguishable. This confirms that the finite-pool noise in the
bootstrap is the dominant source of discrepancy and validates the adequacy of the
Gaussian power approximation for the aggregate test.
Second, the fitted curve, which imposes the
linear-in-$\rho$ structure of Proposition~\ref{prop:heston}, tracks the
simulated power closely across the entire $\rho$ range, demonstrating that for large and moderate leverage effect ($|\rho|\ge0.4$), the rejection probability is close to unity and stable, while for small leverage ($|\rho|<0.4$), power drops and is subject to the data size $B$. Last but not least, the power curve is symmetric in $\rho$, illustrating that the test detects positive and negative leverage with equal probability (two-sided).

\begin{figure}[htbp]
\centering
\includegraphics[width=0.85\textwidth]{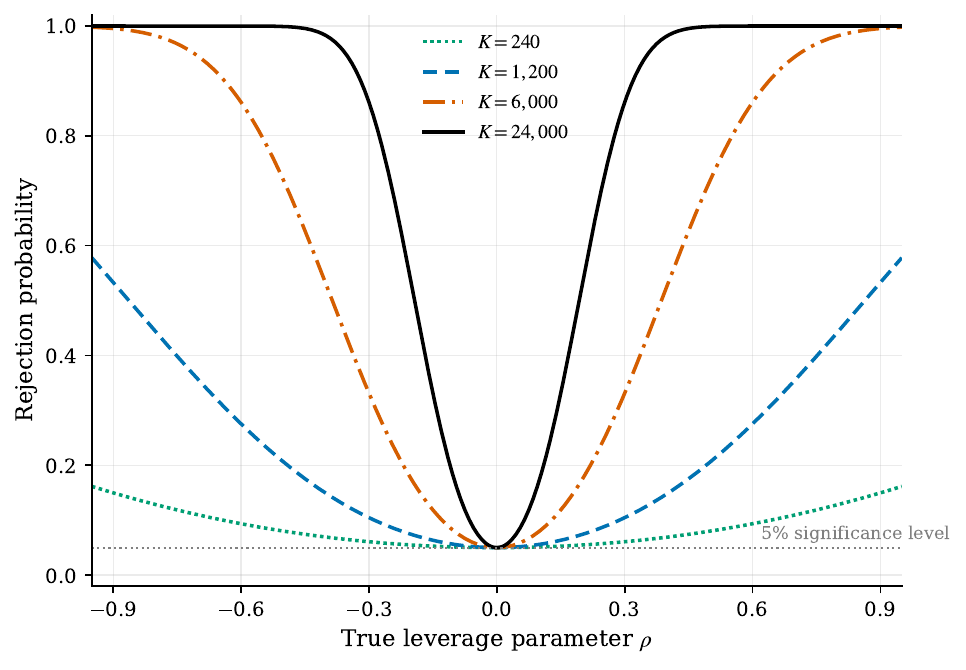}
\captionsetup{justification=raggedright,singlelinecheck=false}
\caption{Analytic power curves of the aggregate mean test over the
leverage parameter $\rho$, w.r.t increasing resample sizes $K$, and fixing the size of the simulation pool at $B=128{,}000$.}
\label{fig:power_vs_K}
\end{figure}

Figure~\ref{fig:power_vs_K} complements the power analysis by fixing
the simulation pool size at $B=128{,}000$ and varying the resampling size
$K\in\{240,\;1{,}200,\;6{,}000,\;24{,}000\}$. Each curve is the analytic
function $\mathrm{Power}(\rho;\,K)$ evaluated at the regression-fitted slope
$\hat a$. As $K$ grows the power envelope tightens around the null:
at $K=240$ power remains below $20\%$ even at $|\rho|=0.9$, while at $K=24{,}000$ the test rejects with near certainty for $|\rho|\ge0.4$.
The transition from low to high power is governed by the product
$|\hat a\,\rho|\sqrt{K}$: the curve reaches $50\%$ power when
$|\hat a\,\rho|\sqrt{K}\approx1.96$. This
confirms that aggregation across windows is the operative lever for
detecting leverage at realistic calibrations.

\begin{remark}[Interpretation of the aggregation $K$]
\label{rem:aggregationinterpretation}
The large values of $K$ required for high power in
Table~\ref{tab:aggpower} and Figure~\ref{fig:power_vs_K} should not be
interpreted as a universal minimum sample size for the proposed test.
Power is governed by the noncentrality
\[
    |\delta|\sqrt{K},
\]
where $\delta=\E[\tau^{0,(1)}]$ is the standardized per-window location
shift under the alternative. Hence, the amount of historical information
required for rejection depends directly on the strength and persistence
of the return asymmetry. Weak alternatives necessarily require long
histories, whereas stronger and more persistent asymmetry can be detected
with substantially smaller $K$.

This feature reflects the hypothesis being tested. Increasing the
within-window sampling frequency makes the realized quadratic
covariation increasingly precise, but does not generate additional
independent realizations of the return-generating mechanism. A single
window therefore contains only limited information about whether the
mechanism itself is systematically asymmetric. By contrast, a
pathwise test based on the infill estimation error can reject with high
probability from a single finely sampled window, but, as discussed in
Section~\ref{sec:twonulls}, such a rejection concerns the realized path
quantity rather than the mechanism-level null considered here. 
For the proposed statistic, power is governed by the standardized mean shift and the number of available windows.

Accordingly, the empirical role of the procedure is not to guarantee
rejection for every asset over a conventional sample period, but to
identify cases in which directional return--variance dependence is
sufficiently strong and persistent to be distinguished from intrinsic
path variation. Failure to reject should therefore not be interpreted as
evidence of symmetry. 
Conversely, rejection provides evidence against the independence null in the direction captured by the aggregate statistic.
Section~\ref{sec:realworld} illustrates this distinction:
with only $K=307$ monthly windows per stock over 1998--2025, the test still
rejects the null for a non-negligible subset of the cross section,
demonstrating that directional asymmetry can be detectable in
realistic samples even though weak alternatives require much larger $K$
for high power.
\end{remark}

\subsection{Bates Framework (jump-present)}\label{sec:bates}
In this section we investigate the impact on the statistics of our jump filter. The jump experiments use the same baseline diffusion calibration
($\kappa=1.5$, $\theta=0.04$, $\xi=0.3$, $\mu=\theta/2$) and overlay
jumps with $\lambda=1.5$, $\mu_J=-0.11$, $\sigma_J=0.12$, the
upper-intensity, negatively-skewed corner of the parameter bounds from
\citet{eraker2004} and hence the most adversarial calibrated case
for the jump filter. Before computing $\QCd_T$, returns flagged by the
\citet{LeeMyklandJump} test (trailing bipower window $K^{LM}=270$, significance $\alpha^{LM}=0.01$) are set to zero, to approximate the continuous component entering the $\mathcal{SVSM}^c$ inference framework of Section~\ref{sec:svsm}. 
The experiments run six configurations,
$32{,}000$ replications each for G--L as shown in
Table~\ref{tab:jumpexps}.

\begin{table}[htbp]
\captionsetup{justification=raggedright,singlelinecheck=false}
\caption{The jump experiment settings.}
\centering
\label{tab:jumpexps}
\begin{tabular}{lll}
\toprule
Config & Other parameters & Verification target \\
\midrule
G\; $\rho=0$    & Baseline jumps, LM-cleaned & Size \\
H\; $\rho=0$    & Baseline jumps, \emph{not} cleaned
  & Jump contamination \\
I\; $\rho=-0.3$ & Baseline jumps, LM-cleaned & Power with jumps removed \\
J\; $\rho=-0.7$ & Baseline jumps, LM-cleaned & Power with jumps removed \\
K\; $\rho=-0.9$ & Baseline jumps, LM-cleaned & Power with jumps removed \\
L\; $\rho=0$    & $\theta\sim\mathcal{U}(0.05,0.5)$,
  $\xi\sim\mathcal{U}(0.1,1)$,
  $\gamma\sim 1+\mathrm{Exp}(4)$, & Size uniformity over \\
  & $\kappa=\gamma\xi^2/(2\theta)$, $\mu=\theta/2$, 
  $\lambda\sim\mathcal{U}(0.25,1.50)$, & parameter uncertainty \\
  & $\mu_J\sim\mathcal{U}(-0.11,0.01)$,
  $\sigma_J\sim\mathcal{U}(0.04,0.12)$, LM-cleaned & \\
\bottomrule
\end{tabular}
\end{table}
Configuration~H serves as a diagnostic: it measures how far the test is displaced when
the jumps are left in, since asymmetric jumps ($\mu_J<0$) generate negative
skewness in the absence of any leverage. Comparing G
with~A isolates the cost of jump removal on size, while comparing I--K with
B--D isolates its cost on power. Configuration~L repeats the parameter-uncertainty
exercise of configuration~F with the jump parameters also drawn at random, over
bounds calibrated to \citet{eraker2004}. The Heston closed form \eqref{eq:qcovdm} provides the continuous martingale-part benchmark for the filtered estimates.


\subsubsection{Size under jumps}
\label{sec:batessize}
\begin{table}[htbp]
\centering
\caption{Jump-present simulation results. Rejections (our method vs.\ BNS method
\citep{Barndorff-N&Shephard2004a}) at two-sided $5\%$ significance.
$t$ denotes the conditional-Gaussian test statistic $t^{0}_n$, whereas
$t_{\mathrm{BNS}}$ denotes the BNS test statistic from \eqref{eq:bnst}.
Column ``flag.'' reports the average number of returns per window flagged by
the \citet{LeeMyklandJump} test as jumps. The last two columns compare the Monte
Carlo mean of the demeaned realized quadratic covariation with its infill
limit in analytical form from Proposition~\ref{prop:heston}, which remains
the benchmark because the filtered path lies in $\mathcal{SVSM}^c$.}
\label{tab:batessize}
\small
\begin{tabular}{lccccccc}
\toprule
Config & mean$(t)$ & sd$(t)$ & rej.\ $t$ & rej.\ BNS
& flag. & $\E[\QCd]$ MC & $\E[L^{0}]$ analytic \\
\midrule
G\; $\rho=0$, jump cleaned & $+0.005$ & $0.992$ & $0.047$        & $0.954$        & $0.133$ & $+3.41\text{e-}7$ & $0$ \\
H\; $\rho=0$, not cleaned  & $-0.071$ & $1.467$ & $0.061$        & $0.864$        & --- & $-5.87\text{e-}5$ & $0$ \\
I\; $\rho=-0.3$            & $-0.017$ & $0.996$ & $0.050^{\ast}$ & $0.054^{\ast}$ & $0.137$ & $-7.48\text{e-}6$ & $-8.08\text{e-}6$ \\
J\; $\rho=-0.7$            & $-0.046$ & $0.998$ & $0.052^{\ast}$ & $0.058^{\ast}$ & $0.136$ & $-1.85\text{e-}5$ & $-1.89\text{e-}5$ \\
K\; $\rho=-0.9$            & $-0.066$ & $0.997$ & $0.053^{\ast}$ & $0.065^{\ast}$ & $0.135$ & $-2.45\text{e-}5$ & $-2.42\text{e-}5$ \\
L\; random parameters      & $-0.004$ & $1.000$ & $0.052$        & $0.953$        & $0.083$ & --- & --- \\
\bottomrule
\end{tabular}

\vspace{2pt}
{\footnotesize\raggedright
$^{\ast}$ Size-adjusted rejection rates. For the power configurations
(I--K) the entries in columns rej.\ $t$ and rej.\ BNS use size-adjusted
critical values instead of $\pm1.96$: the empirical $2.5$th and $97.5$th
percentiles of each $t$-statistic under the size configuration~G ($\rho=0$,
jumps present and cleaned) are taken as critical values, and the rejection
rate is the fraction of the $t$-statistic falling outside them. Adjusting
against~G rather than against the jump-free configuration~A ensures that
any residual distortion contributed by the filter itself is absorbed into
the reference distribution rather than counted as power. Configurations
G, H, and L have $\rho=0$ and retain the $5\%$ $\pm1.96$ critical
values.\par}
\end{table}

Table~\ref{tab:batessize} reports the jump-present simulation results. First, the
conditional-Gaussian statistic retains the correct size once jumps are
removed: rejection rates of $4.7\%$ and $5.2\%$ in the null configurations
G and~L, with mean close to zero and standard deviation close to one, are
indistinguishable from the jump-free configurations A and~F of
Table~\ref{tab:size}. The detector flags approximately $0.13$ returns per window, close to the expected jump count $\lambda T=0.125$, while the resulting test retains near-nominal size in these experiments.

Second, configuration~H shows the influence of jumps on the tests, and
the damage is larger than its rejection rate suggests. Although the raw rejection rate remains modest at $6.1\%$, the jumps shift
the mean of $t^{0}_n$ to $-0.071$, almost the same size as the shift induced by
true leverage at $\rho=-0.9$ (configuration~K, $-0.066$), and
displace $\E[\QCd_T]$ to $-5.87\text{e-}5$ against a true value of zero,
roughly two times the magnitude of the $\rho=-0.9$ signal itself. The
dispersion simultaneously inflates, $\mathrm{sd}(t^{0}_n)$ rising from
$0.992$ to $1.467$. The modest rejection rate is thus misleading: the
contamination is concentrated in rare and large realizations, so the bulk of
the distribution remains close to $\mathcal{N}(0,1)$ while its first moment
is badly displaced. Since the aggregate mean test of Section~\ref{sec:power}
targets exactly that first moment, uncleaned jumps would be amplified by
$\sqrt{K}$ into spurious rejections of $H_0$ and read as persistent asymmetry, though the true leverage is zero.  

Third, the Monte Carlo means of $\QCd_T$ reproduce the closed form
\eqref{eq:qcovdm} across all leverage strengths, with deviations of the same order as in the jump-free configurations B--D.
The filtered estimates remain close to the continuous benchmark, with discrepancies comparable to those in the jump-free experiments.

Last but not least, under the alternative (I--K) neither statistic has useful per-window
power, with size-adjusted rejection rates between $5.0\%$ and $6.5\%$,
exactly as in the jump-free configurations B--D. The
aggregate mean test of Section~\ref{sec:power} remains the route to power,
and we verify in Section~\ref{sec:batespower} that its properties carry
over.

\subsubsection{Power and jump intensity}
\label{sec:batespower}
We repeat the aggregation exercise of Section~\ref{sec:power} in this new context. The construction is identical: $K$ statistics are drawn with
replacement from the simulated pool and we reject when
$|\sqrt K\,\bar t^{\,0}_{K,n}|>z_{0.975}$, except for one change. The recentering uses
the null pool of configuration~G rather than that of configuration~A, so
that any residual displacement contributed by the \citet{LeeMyklandJump}
filter is absorbed into the null instead of being counted as power. 

\begin{table}[htbp]
\centering
\caption{Power of the aggregate mean test $\sqrt K\,\bar t^{\,0}_{K,n}$
(conditional-Gaussian statistic) under the jump pipeline, at the two-sided
$5\%$ significance level. Each entry is the bootstrap rejection rate
over $16{,}000$ draws from a pool of $32{,}000$ simulated windows. Config~G
($\rho=0$, jumps cleaned) is the size row and supplies the recentering.}
\label{tab:batesaggpower}
\small
\begin{tabular}{lcccccc}
\toprule
Config & $K=60$ & $K=120$ & $K=240$ & $K=1{,}200$ & $K=6{,}000$ & $K=24{,}000$ \\
\midrule
G\; $\rho=0$, jump cleaned     & $0.049$ & $0.048$ & $0.048$ & $0.049$ & $0.047$ & $0.050$ \\
I\; $\rho=-0.3$           & $0.051$ & $0.059$ & $0.063$ & $0.118$ & $0.404$ & $0.927$ \\
J\; $\rho=-0.7$           & $0.065$ & $0.084$ & $0.125$ & $0.426$ & $0.977$ & $1.000$ \\
K\; $\rho=-0.9$           & $0.085$ & $0.117$ & $0.193$ & $0.690$ & $1.000$ & $1.000$ \\
\midrule
H\; $\rho=0$, not cleaned & $0.165$ & $0.197$ & $0.270$ & $0.673$ & $0.998$ & $1.000$ \\
\bottomrule
\end{tabular}
\end{table}

Table~\ref{tab:batesaggpower} carries three messages. First, configuration~G likewise gives near-nominal null rejection rates under the recentered resampling scheme with jump filtering. 

Second, configuration~H quantifies the impact of jump contamination by showing what the same aggregation does to uncleaned
data, and the contrast is obvious. 
%
Although the continuous leverage parameter is zero, unremoved asymmetric jumps violate the maintained continuous-price assumption; rejection rates reach $67.3\%$ at $K=1{,}200$ and $99.8\%$ at $K=6{,}000$.

The
jump-induced displacement of $\E[t^{0}_n]$ documented in
Section~\ref{sec:batessize} is amplified by $\sqrt{K}$. As a result, testing
for persistent asymmetry by aggregation while omitting jump removal would not only weaken the test but also invert its conclusion, reporting
leverage with near-certainty where the data-generating process has none.

Third, the rejection rate rises monotonically in $K$ and in $|\rho|$ at every leverage strength considered, reaching $92.7\%$, $100\%$, and $100\%$ at $K=24{,}000$ for $\rho=-0.3$, $-0.7$, and $-0.9$, closely matching the jump-free rates of Table~\ref{tab:aggpower}. Across the investigated calibrations, jump removal using the
\citet{LeeMyklandJump} method yields power comparable to the jump-free benchmark
and preserves the gains from aggregation documented in Section~\ref{sec:power}. 

\begin{figure}[H]
\centering
\includegraphics[width=0.85\textwidth]{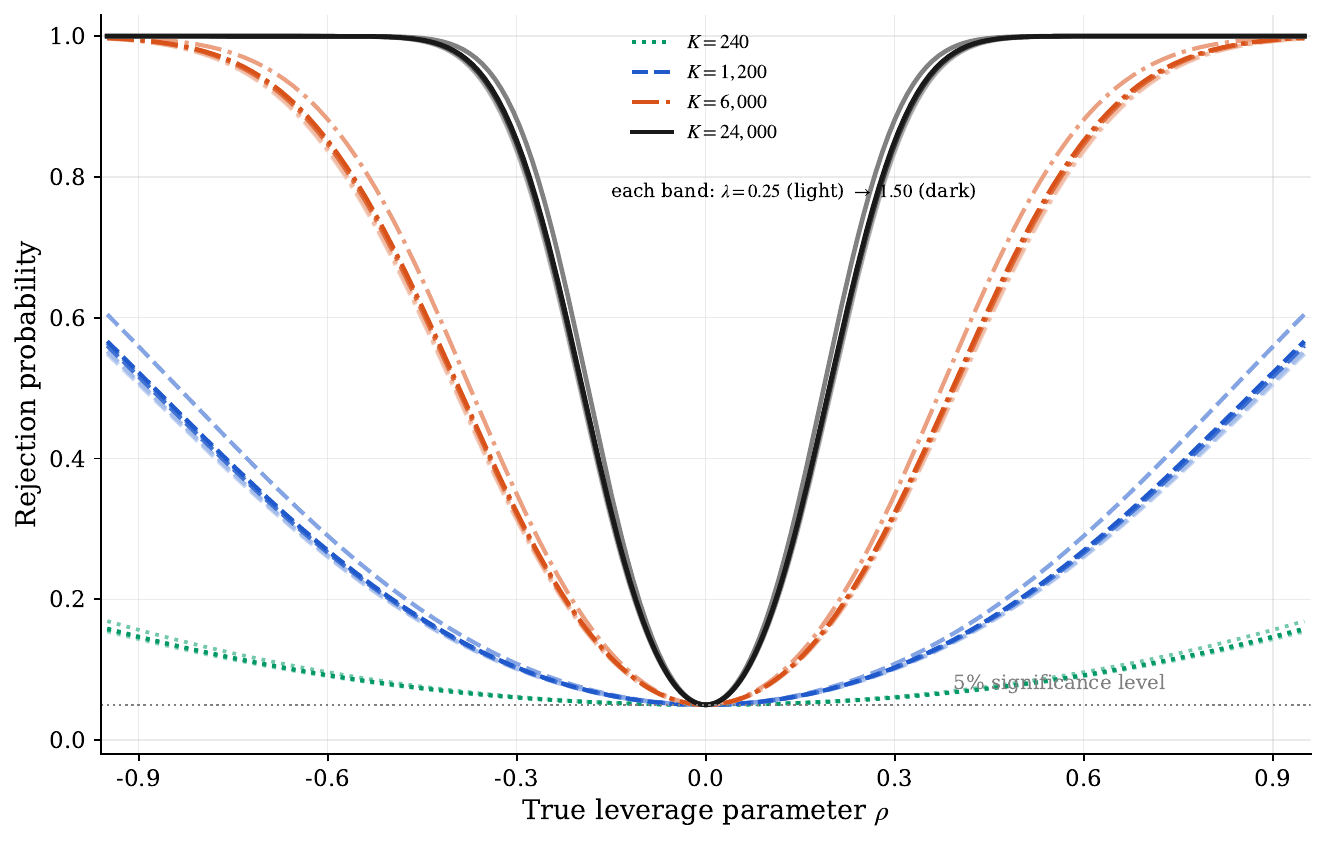}
\caption{Analytic power curves of the aggregate mean test over the leverage
parameter $\rho$ for four aggregate panel sizes
$K\in\{240,\;1{,}200,\;6{,}000,\;24{,}000\}$, each evaluated under six
jump intensities $\lambda\in\{0.25,0.50,0.75,1.00,1.25,1.50\}$ with darkness
increasing in $\lambda$. The per-window slope $\hat a(\lambda)$ is fitted
from $128{,}000$ simulated windows per $(\lambda,\rho)$.}
\label{fig:jump_power_4K}
\end{figure}

Figure~\ref{fig:jump_power_4K} extends this comparison to the full power
surface. We run the same regression analysis as in Section~\ref{sec:power} but additionally with respect to six jump intensity levels, i.e., $\lambda\in\{0.25,0.50,0.75,1.00,1.25,1.50\}$. As a result, this yields a fitted analytic power curve at aggregate size $K$ and jump intensity $\lambda$:
\[
\mathrm{Power}(\rho;\,K,\lambda)
\;=\;\Phi\bigl(-z_{0.975}+|\hat a(\lambda)\,\rho|\sqrt{K}\bigr)
\;+\;\Phi\bigl(-z_{0.975}-|\hat a(\lambda)\,\rho|\sqrt{K}\bigr),
\]
where $\hat a(\lambda)$ is estimated from $128{,}000$
windows per $(\lambda,\rho)$ on a $19$-point $\rho$ grid. Within each $K$
band ($K\in\{240,\;1{,}200,\;6{,}000,\;24{,}000\}$), the six curves, corresponding to jump intensities $\lambda=0.25,\,0.50,\,\ldots,\,1.50$, are visually indistinguishable: the
fitted slopes range from $0.0633$ to $0.0676$, all within roughly two standard errors
of the jump-free slope estimate $\hat a=0.0655$ ($\mathrm{se}\;0.001$) of
Section~\ref{sec:power}, and a trend regression of $\hat a$ on $\lambda$
yields a slope not significantly different from zero. The fitted power curves show little sensitivity to jump intensity over the calibrated range.


\section{Real-world Application}\label{sec:realworld}

\begin{figure}[H]
\centering
\includegraphics[width=\textwidth]{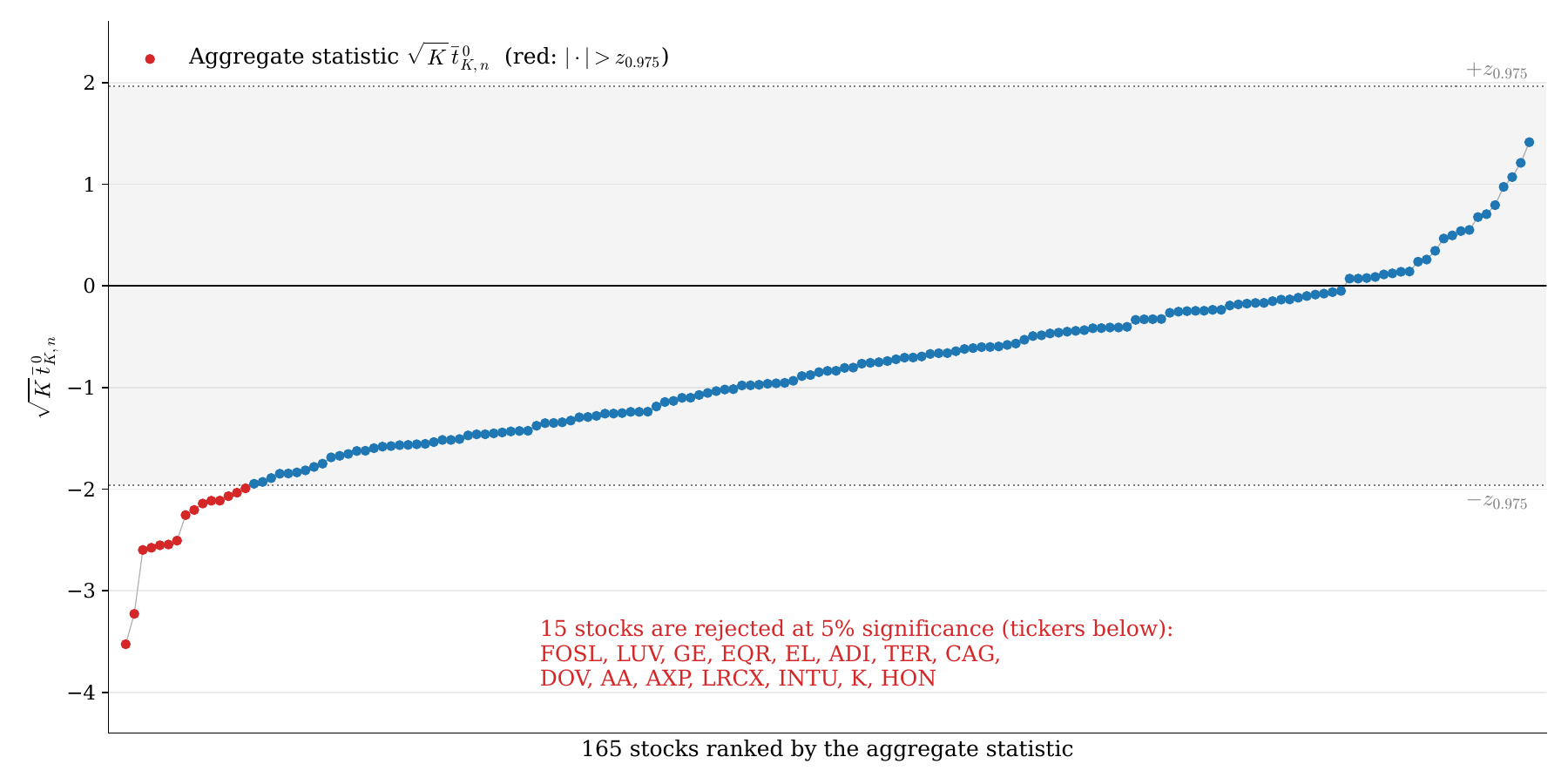}
\captionsetup{justification=raggedright,singlelinecheck=false}
\caption{Aggregate mean test $\sqrt K\,\bar t^{\,0}_{K,n}$ ($K=307$) for $165$
U.S.~500 stocks, ranked from most negative to most positive. Filled
circles are the conditional-Gaussian statistic, red where
$|\sqrt K\,\bar t^{\,0}_{K,n}|>z_{0.975}$. The 15 rejected tickers
are listed. The shaded band is the two-sided $5\%$ acceptance region.}
\label{fig:rw_ranked}
\end{figure}

In this section, we apply the proposed hypothesis test to U.S. stock market data, sourced from the TAQ database. Our sample includes $165$ highly liquid stocks with a complete record of $5$-minute intraday returns from
April~1998 to April~2025.\footnote{From the S\&P 500, we retain stocks with five-minute return data from 1998 to 2025, resulting in 165 stocks. We remove overnight returns and exclude data from early-closing days.} For each stock the statistic is formed exactly
as in the simulation pipeline: returns are grouped into non-overlapping
monthly windows ($22$ trading days, $78$ intraday returns per day, so $n=1{,}716$), jumps are detected and removed
by the \citet{LeeMyklandJump} test, returns are demeaned within the
window, and the bridge-corrected conditional-Gaussian statistic
$t^{0}_{n}$ of Theorem~\ref{thm:tstat} is computed. This produces a panel of $K=307$ monthly statistics per stock. Applying the aggregate mean test
of Corollary~\ref{cor:powerdm} stock by stock,
Figure~\ref{fig:rw_ranked} ranks the $165$ realizations of
$\sqrt K\,\bar t^{\,0}_{K,n}$: $H_0$ is rejected for 15 stocks ($9.1\%$) at the $5\%$ significance level, and every one lies in the
negative tail. None of them shows a significantly positive asymmetry. The reported significance levels are stock-specific and are not adjusted for multiple testing.

\section{Conclusion} \label{sec:conclusion}
To conclude, this paper develops inference for the realized quadratic covariation between
a return path and its running square, $[R,R^{2}]$, as a measure of
return asymmetry under the continuous stochastic volatility semimartingale
class $\mathcal{SVSM}^c$. \citet{Barndorff-N&Shephard2004a} were the first to provide an asymptotic distribution theory for the realized covariation matrix, studentized by the variance of the infill estimation error, which supports confidence intervals for the realized quadratic covariation and hypothesis tests of whether the realized statistic has converged to its probability limit.

We complement their distribution theory by considering the conditional variance of the limit itself given the volatility path. Under the volatility--Brownian independence null and the maintained drift assumptions, the limiting quadratic covariation divided by its conditional standard deviation is exactly standard normal, yielding a feasible asymptotically pivotal test. 
Moreover, we show that an empirically demeaned version of the test statistic, which removes an unknown but window-constant drift, preserves asymptotic normality with a bridge correction. Extending the
analysis across non-overlapping windows yields the power theory of the
paper. 

Under the null, the limiting standardized window statistics are independent standard normals, and the feasible aggregate mean test has asymptotic size $\alpha$ for every fixed $K$.
%
%
Under stationary
and ergodic window blocks, whenever the directional asymmetry produces a
nonzero population shift
$\delta=\E[\tau^{0,(1)}]$, the sample mean converges to $\delta$ and the
test is consistent as first $n\to\infty$ and then $K\to\infty$. 
Thus, power accumulates through the
noncentrality $|\delta|\sqrt{K}$: weak asymmetry requires a long history,
whereas stronger and more persistent asymmetry is detectable with a
smaller number of windows. This is a consequence of testing a
mechanism-level hypothesis rather than a realized-path hypothesis.

We also validate the size and the power of the proposed test using a simulation study concerning both jump-free and jump-present scenarios. In addition, we apply this test to $165$ highly liquid U.S.\ stocks over 1998--2025.  Using $K=307$ non-overlapping monthly windows per stock, the aggregate mean test rejects the null for $15$ stocks, all in the
negative direction, while the remaining stocks do not reject the independence null. 
A failure to reject
should not be interpreted as evidence that return asymmetry is absent.
Rather, rejection identifies cases in which directional asymmetry is
sufficiently strong and persistent to dominate the intrinsic variation
of the realized path.

The proposed aggregate mean test is intentionally targeted rather than
omnibus. It is consistent against alternatives for which
$\delta\neq0$, that is, directional return--variance dependence producing
a nonzero signed shift in the standardized quadratic covariation. More
general forms of dependence between returns and volatility innovations may violate
the independence null while leaving $\delta=0$ through offsetting or
sign-symmetric effects. Extending the present framework to tests that
exploit higher-order features of the window-statistic distribution, for
example smooth tests based on higher Hermite coordinates, provides a
natural direction for future research. A second potential extension is to develop
joint cross-sectional inference that exploits common asymmetry across
assets while explicitly accounting for market-wide dependence.

\newpage
\bibliographystyle{elsarticle-harv} 
{\small
\renewcommand{\baselinestretch}{1} 
\bibliography{literature}
}

\renewcommand{\baselinestretch}{1.5} 
\newpage
\appendix
\renewcommand{\thesubsection}{\Alph{section}.\roman{subsection}}

\section{Appendix}\label{Appendix}
\subsection{Proofs of the main results}\label{app:proofs}
This appendix collects the proofs of the results stated in
Sections~\ref{sec:test} and~\ref{sec:empimpl}, in the order in which the
results appear. The proofs of Lemmas~\ref{lem:limitraw}
and~\ref{lem:driftremoval} are short and remain in the main text.

\subsubsection{Proofs for Section~\ref{sec:mart}: the martingale case}\label{app:proofs-mart}

\phantomsection\label{prf:riemann}
\begin{proof}[Proof of Lemma~\ref{lem:riemann}]
Work conditionally on $\Gsig$. Decompose the increment into its drift and
martingale parts, $\Delta R_j=\Delta A_j+\Delta M_j$, where
$\Delta A_j=\int_{t_{j-1}}^{t_j}\alpha_u\,du$ and
$\Delta M_j=\int_{t_{j-1}}^{t_j}\sigma_u\,dW_u$. Under $H_0$,
$\sigma\perp W$, so given $\F_{j-1}\vee\Gsig$ the martingale increment
$\Delta M_j$ is Gaussian with mean zero and variance
$\Delta\Lam_j=\Lam_j-\Lam_{j-1}$. Write
$(\Delta R_j)^2=\Delta\Lam_j+\eta_j+\zeta_j$, where
$\eta_j=\Delta M_j^2-\Delta\Lam_j$ collects the centered-square
fluctuation, with $\E[\eta_j\mid\F_{j-1}\vee\Gsig]=0$ and
$\E[\eta_j^2\mid\F_{j-1}\vee\Gsig]=2\Delta\Lam_j^2$, and
$\zeta_j=2\Delta A_j\,\Delta M_j+\Delta A_j^2$ collects the drift
contribution. Then
\[
\sum_j f_{j-1}(\Delta R_j)^2
=\sum_j f_{j-1}\Delta\Lam_j+\sum_j f_{j-1}\eta_j+\sum_j f_{j-1}\zeta_j .
\]
The first sum converges a.s.\ to $\int_0^T f_u\,d\Lam_u$ as a
Riemann--Stieltjes sum, $f$ being continuous and $\Lam$ continuous and
increasing. The second is a martingale-difference sum with
\[
\E\Bigl[\Bigl(\textstyle\sum_j f_{j-1}\eta_j\Bigr)^{\!2}\Bigm|\Gsig\Bigr]
=2\sum_j \E\bigl[f_{j-1}^2\mid\Gsig\bigr]\Delta\Lam_j^2
\;\le\; 2\,\E\Bigl[\sup_u f_u^2\Bigm|\Gsig\Bigr]\,
\Bigl(\max_j\Delta\Lam_j\Bigr)\Lam_T \;\longrightarrow\;0 ,
\]
because $\Lam$ is continuous, so $\max_j\Delta\Lam_j\to0$.

For the third sum no martingale structure is required: it vanishes by
pathwise bounds alone, so no condition relating $\alpha$ to $W$ is
imposed. Since $\alpha$ is c\`adl\`ag on $[0,T]$, it is bounded, say
$|\alpha_u|\le C$, whence $|\Delta A_j|\le C\Delta_n$ pathwise. The
drift-squared term obeys
\[
\Bigl|\sum_j f_{j-1}\,\Delta A_j^2\Bigr|
\le \Bigl(\max_j\Delta A_j^2\Bigr)\sum_j|f_{j-1}|
\le C^2\Delta_n^2\,n\,\sup_{u\le T}|f_u|
= C^2 T\,\Delta_n\,\sup_{u\le T}|f_u| \;\longrightarrow\;0
\]
a.s., since $\sup_{u\le T}|f_u|<\infty$ a.s.\ ($f$ continuous on a
compact). For the cross term, bound
\[
\Bigl|2\sum_j f_{j-1}\,\Delta A_j\,\Delta M_j\Bigr|
\;\le\; 2\,C\Delta_n\,\sup_{u\le T}|f_u|\;\sum_j|\Delta M_j| ,
\]
and control the last factor by the Cauchy--Schwarz inequality and the
conditional It\^o isometry:
\[
\E\Bigl[\sum_j|\Delta M_j|\Bigm|\Gsig\Bigr]
\;\le\;\sqrt{n}\,\Bigl(\sum_j\E\bigl[\Delta M_j^2\mid\Gsig\bigr]\Bigr)^{\!1/2}
=\sqrt{n\,\Lam_T}
=\sqrt{T\,\Lam_T}\;\Delta_n^{-1/2} ,
\]
so that, by Markov's inequality applied conditionally on $\Gsig$,
$\sum_j|\Delta M_j|=O_p(\Delta_n^{-1/2})$ and hence
$\Delta_n\sum_j|\Delta M_j|=O_p(\Delta_n^{1/2})\plim0$. Since
$\sup_{u\le T}|f_u|<\infty$ a.s., the cross term satisfies
\[
\Bigl|2\sum_j f_{j-1}\,\Delta A_j\,\Delta M_j\Bigr|
\;\le\; 2C\,\Delta_n\,\sup_{u\le T}|f_u|\,\sum_j|\Delta M_j|
\;\plim\;0 .
\] 
Conditional convergence in
probability for a.e.\ volatility path implies unconditional convergence
in probability.
\end{proof}



\phantomsection\label{prf:righttag}
\begin{proof}[Proof of Corollary~\ref{cor:righttag}]
Since $f$ is continuous on $[0,T]$, hence uniformly continuous,
$\max_j|f_j-f_{j-1}|\to0$; with $\sum_j(\Delta R_j)^2\plim\Lam_T$
(Lemma~\ref{lem:riemann}, $f\equiv\mathbf 1$),
\[
\Bigl|\sum_j (f_j-f_{j-1})(\Delta R_j)^2\Bigr|
\le \max_j|f_j-f_{j-1}|\sum_j(\Delta R_j)^2 \plim 0 ,
\]
so the left- and right-endpoint sums share the limit of Lemma~\ref{lem:riemann}.
\end{proof}

\phantomsection\label{prf:condgaussraw}
\begin{proof}[Proof of Proposition~\ref{prop:condgaussraw}]
By Fubini's theorem for stochastic integrals (justified pathwise,
conditionally on $\Gsig$, since $R$ is then a Gaussian martingale and
$\Lam$ a deterministic finite measure),
$\int_0^T R_u\,d\Lam_u=\int_0^T\!\int_0^u dR_s\,d\Lam_u
=\int_0^T(\Lam_T-\Lam_s)\,dR_s$. Conditionally on $\Gsig$, $h$ is a
deterministic bounded function and $R$ a continuous Gaussian martingale
with variance function $\Lam$. The Wiener integral $2\int h\,dR$ is
therefore Gaussian with mean zero and variance $4\int h^{2}\,d\Lam$.
Positivity holds because $h(u)>0$ for $u<T$ and $\Lam$ is strictly
increasing under Assumption~\ref{ass:sv}.
\end{proof}

\phantomsection\label{prf:mart}
\begin{proof}[Proof of Theorem~\ref{thm:mart}]
All convergence arguments below are conducted under the regular conditional law given $\mathcal G$, for almost every volatility path.

\emph{(i) Uniform consistency of $\widehat\Lam$.} By
Lemma~\ref{lem:riemann} ($f\equiv\mathbf 1$ on $[0,t]$),
$\widehat\Lam_n(t)\plim\Lam_t$ pointwise. By
Assumption~\ref{ass:sv}, $\Lam$ is Lipschitz with constant $\overline\sigma^2$, i.e., $\Lam_u-\Lam_s\le\overline\sigma^{2}(u-s)$ for $0\le s<u\le T$. Then for
any $\varepsilon>0$ a deterministic grid $0=t_0<\dots<t_m=T$ can be chosen with
$\Lam_{t_i}-\Lam_{t_{i-1}}<\varepsilon$ on every path. Since $\widehat\Lam_n$
and $\Lam$ both are nondecreasing, the usual finite-grid argument gives
\[
\sup_{t\le T}\bigl|\widehat\Lam_n(t)-\Lam_t\bigr|
\le\max_{i\le m}\bigl|\widehat\Lam_n(t_i)-\Lam_{t_i}\bigr|+\varepsilon ,
\]
and a union bound over the $m$ grid points yields
$\sup_{t\le T}|\widehat\Lam_n(t)-\Lam_t|\plim0$.


\emph{(ii) Consistency of $\widehat V$.} Let $g(u)=(\Lam_T-\Lam_u)^{2}$ and
$\widehat g_n(u)=(\widehat\Lam_n(T)-\widehat\Lam_n(u))^{2}$; by
(i) and continuous mapping, $\sup_u|\widehat g_n(u)-g(u)|\plim0$. Then
\[
\tfrac14\widehat V_n
=\underbrace{\sum_j g(t_j)(\Delta R_j)^2}_{\to\,\int_0^T g\,d\Lam
\ \text{(Corollary \ref{cor:righttag})}}
+\underbrace{\sum_j\{\widehat g_n(t_j)-g(t_j)\}(\Delta R_j)^2}_{
\le\,\sup_u|\widehat g_n-g|\cdot\widehat\Lam_n(T)\,\plim\,0} .
\]
Hence $\widehat V_n\plim V$.

\emph{(iii) The $t$-statistic.} By Lemma~\ref{lem:limitraw},
$\QC_T=L+o_p(1)$. By (i)--(ii) and $V>0$ a.s.,
$(\widehat V_n)^{-1/2}=(V)^{-1/2}+o_p(1)$. Slutsky gives
$t_n=L/\sqrt{V}+o_p(1)$, and by
Proposition~\ref{prop:condgaussraw} the ratio $L/\sqrt{V}$ is
\emph{exactly} standard normal conditionally on $\Gsig$, i.e., $V$ is
$\Gsig$-measurable, so the self-normalization removes the volatility path
entirely. Hence $\Prob(t_n\le z\mid\Gsig)\to\Phi(z)$ for a.e.\ path,
and unconditionally by dominated convergence.
\end{proof}

\phantomsection\label{prf:condindep}
\begin{proof}[Proof of Lemma~\ref{lem:condindep}]
Let $Y_k:=(W_u-W_{(k-1)T})_{u\in\mathcal I_k}$ and
$S_k:=(\sigma_u)_{u\in\mathcal I_k}$, $k=1,\dots,K$. Under
Assumption~\ref{ass:mart}, $R^{(k)}_u=\int_{(k-1)T}^{u}\sigma_s\,dW_s$, so
there are measurable maps $\Phi_{k,n}$ and $\Psi_k$ with
\[
\bigl(\QC^{(k)},\widehat V^{(k)}_n\bigr)=\Phi_{k,n}(S_k,Y_k),\qquad
\bigl(L^{(k)},V^{(k)}\bigr)=\Psi_k(S_k,Y_k),
\]
where the second map is meant in the sense of the stochastic integral
$L^{(k)}=2\int_{\mathcal I_k}h^{(k)}\,dR^{(k)}$ with
$h^{(k)}(u)=\Lam_{kT}-\Lam_u$, and $V^{(k)}=4\int_{\mathcal I_k}(h^{(k)})^2d\Lam$
depends on $S_k$ only. Each $S_k$ is $\Gsig_K$-measurable. Under $H_0$
the whole path $W$ is independent of $\Gsig_K$, hence so is
$(Y_1,\dots,Y_K)$. And $Y_1,\dots,Y_K$ are mutually independent, being
increments of a Brownian motion over disjoint intervals. Therefore the
conditional law of $(Y_1,\dots,Y_K)$ given $\Gsig_K$ is its unconditional
law, the product measure $\bigotimes_{k=1}^{K}P_{Y_k}$, where $P_{Y_k}$
denotes the law of $Y_k$ and $\bigotimes$ denotes the product measure under which
the coordinates are independent with the prescribed marginals. By the freezing
lemma for conditional expectations, for bounded measurable
$g_1,\dots,g_K$,
\[
\E\Bigl[\prod_{k=1}^{K}g_k\bigl(\Phi_{k,n}(S_k,Y_k)\bigr)\Bigm|\Gsig_K\Bigr]
=\int\prod_{k=1}^{K}g_k\bigl(\Phi_{k,n}(S_k,y_k)\bigr)\,\bigotimes_{k}P_{Y_k}(dy_k)
=\prod_{k=1}^{K}\E\bigl[g_k\bigl(\Phi_{k,n}(S_k,Y_k)\bigr)\bigm|\Gsig_K\bigr]
\]
a.s., where the last step is derived by Fubini. This is the conditional independence of
the pairs $(\QC^{(k)},\widehat V^{(k)}_n)$ across $k$, for every $n$.

The same computation with $\Psi_k$ in place of $\Phi_{k,n}$ shows that
$L^{(1)},\dots,L^{(K)}$ are conditionally independent given $\Gsig_K$, and
that the conditional law of $L^{(k)}$ given $\Gsig_K$ is the law of
$2\int_{\mathcal I_k}h^{(k)}\,dR^{(k)}$ with $S_k$ frozen at its realized
value. With $S_k$ frozen, $h^{(k)}$ and $\sigma$ are deterministic on
$\mathcal I_k$ and $R^{(k)}$ is a Wiener integral, so this law is
$\mathcal N(0,V^{(k)})$. Since $V^{(k)}$ is $\Gsig_K$-measurable and
positive, $\tau^{(k)}\mid\Gsig_K\sim\mathcal N(0,1)$, and conditional
independence gives $(\tau^{(1)},\dots,\tau^{(K)})\mid\Gsig_K\sim\mathcal N(0,I_K)$.
As this law does not depend on the conditioning variables, it is also the
unconditional law and the vector is independent of $\Gsig_K$.
\end{proof}

\phantomsection\label{prf:power}
\begin{proof}[Proof of Corollary~\ref{cor:power}]
(i) Lemma~\ref{lem:condindep} gives the exact law of $\sqrt K\,\bar\tau_K$.
For fixed $K$, $\sqrt K(\bar t_{K,n}-\bar\tau_K)$ is a finite sum of
$o_p(1)$ terms, and Slutsky's lemma applies.
(ii) Birkhoff's ergodic theorem gives $\bar\tau_K\to\delta$ a.s., so
$|\sqrt K\,\bar\tau_K|\to\infty$ a.s.\ when $\delta\neq0$. For fixed $K$,
$\sqrt K\,\bar t_{K,n}\dlim\sqrt K\,\bar\tau_K$ as $n\to\infty$ by Lemma~\ref{lem:Vhatalt}. Since $\{|x|>z_{1-\alpha/2}\}$ is open, the
Portmanteau theorem gives
$\liminf_n\Prob(|\sqrt K\,\bar t_{K,n}|>z_{1-\alpha/2})
\ge\Prob(|\sqrt K\,\bar\tau_K|>z_{1-\alpha/2})$. Hence
$\lim_{K\to\infty}\liminf_{n\to\infty}
\Prob(|\sqrt K\,\bar t_{K,n}|>z_{1-\alpha/2})=1$. 
\end{proof}

\subsubsection{Proofs for Section~\ref{sec:empimpl}: the demeaned statistic}\label{app:proofs-dm}

\phantomsection\label{prf:limitdm}
\begin{proof}[Proof of Lemma~\ref{lem:limitdm}]
For the cubic term,
$\bigl|\sum_j(\Delta\Rt_j)^3\bigr|\le\max_j|\Delta\Rt_j|
\sum_j(\Delta\Rt_j)^2$. Here
\[
    \sum_{j=1}^n(\Delta\widetilde R_j)^2=\sum_{j=1}^n(\Delta R_j)^2-\frac{R_T^2}{n} \xrightarrow{p}\Lambda_T,
\]
by the classical convergence \citep{jacod1987limit} of realized quadratic variation for $R$
%
%
, while $\max_j|\Delta\Rt_j|\plim0$ by the uniform continuity of $\Rt$. 
Hence the cubic sum is $\plim0$. For the leading term, expand
$(\Delta\Rt_j)^2=(\Delta R_j)^2-2(R_T/n)\Delta R_j+R_T^2/n^2$:
\[
\sum_j 2\Rt_{j-1}(\Delta\Rt_j)^2
=\underbrace{\sum_j 2\Rt_{j-1}(\Delta R_j)^2}_{(a)}
-\underbrace{\frac{4R_T}{n}\sum_j \Rt_{j-1}\Delta R_j}_{(b)}
+\underbrace{\frac{2R_T^2}{n^2}\sum_j \Rt_{j-1}}_{(c)} .
\]
Term $(b)$: $\sum_j\Rt_{j-1}\Delta R_j\plim
\int_0^T R_u\,dR_u-(R_T/T)\int_0^T u\,dR_u=O_p(1)$, so $(b)=O_p(n^{-1})$.
Term $(c)\le2R_T^2\sup_u|\Rt_u|/n=O_p(n^{-1})$. For term $(a)$, we make the following decomposition and consider the anticipative component factor including $R_T$ as a
fixed random scalar multiplying an adapted sum,
\[
(a) \;=\; 2\sum_j R_{j-1}(\Delta R_j)^2
\;-\;\frac{2R_T}{T}\sum_j t_{j-1}(\Delta R_j)^2
\;\plim\;2\int_0^T R_u\,d\Lam_u-\frac{2R_T}{T}\int_0^T u\,d\Lam_u ,
\]
by the classical convergence of weighted realized quadratic variation
\citep{jacod1987limit} applied to the adapted, continuous integrand
$f=R$ and the deterministic integrand $f(u)=u$. Collecting terms yields
$2\int_0^T\Rt_u\,d\Lam_u$. Since $u\mapsto(u/T)R_T$ is a continuous
finite-variation path, $[\Rt]=[R]$ and
$[\Rt,\Rt^2]_T=\int_0^T2\Rt_u\,d[\Rt]_u$, which is \eqref{eq:limit}.
\end{proof}

\phantomsection\label{prf:condgaussdm}
\begin{proof}[Proof of Proposition~\ref{prop:condgaussdm}]
By Lemma~\ref{lem:driftremoval} the demeaned returns, and hence $L^{0}$ and
$V^{0}$, are pathwise identical to those computed from the martingale part
$M$. We may therefore take $\alpha\equiv0$, so that $R$ is a continuous
local martingale with $[R]=\Lam$. The Fubini step of
Proposition~\ref{prop:condgaussraw} yields
$\int_0^TR_u\,d\Lam_u=\int_0^T(\Lam_T-\Lam_s)\,dR_s$. The bridging term
contributes $\int_0^T(u/T)R_T\,d\Lam_u=c\,R_T=c\int_0^TdR_s$. Subtracting
yields \eqref{eq:rep}, and conditional Gaussianity with variance
\eqref{eq:V} follows exactly as before. For positivity: under
Assumption~\ref{ass:sv}, $\Lam$ is strictly increasing, so $h^{0}$ is
continuous and strictly decreasing with $h^{0}(0)=\Lam_T-c>0$ (since
$c=T^{-1}\!\int u\,d\Lam<\Lam_T$) and $h^{0}(T)=-c<0$ (since $c>0$).
Hence $h^{0}$ vanishes at exactly one point and $V^{0}=4\int (h^{0})^2\,d\Lam>0$.
\end{proof}

\phantomsection\label{prf:consistency}
\begin{proof}[Proof of Theorem~\ref{thm:consistency}]
By Lemma~\ref{lem:driftremoval} the demeaned returns, and hence
$\widehat V^{0}_n$, are pathwise identical to those computed from the
martingale part. We may therefore take $\alpha\equiv0$.

\emph{Step 1: demeaning is negligible in squared returns.} With
$\bar r=R_T/n$,
\[
\sum_{j}\bigl|\rt_j^{\,2}-(\Delta R_j)^2\bigr|
\le 2|\bar r|\sum_j|\Delta R_j|+n\bar r^{\,2}
= O_p\bigl(n^{-1}\cdot n^{1/2}\bigr)+O_p\bigl(n^{-1}\bigr)
= O_p\bigl(n^{-1/2}\bigr),
\]
since $\sum_j|\Delta R_j|=O_p(\sqrt{n})$ under Assumption~\ref{ass:sv}.
Hence all sums below may be computed with $(\Delta R_j)^2$ in place of
$\rt_j^{\,2}$ at cost $o_p(1)$, uniformly over bounded weights.

\emph{Step 2: ingredients.} Part (i) of the proof of
Theorem~\ref{thm:mart} gives
$\sup_{t\le T}|\widehat\Lam^{0}_n(t)-\Lam_t|\plim0$.
Corollary~\ref{cor:righttag} with $f_u=u/T$ gives $\widehat c\plim c$.

\emph{Step 3: conclusion.} With $g(u)=\{(\Lam_T-\Lam_u)-c\}^2$ and
$\widehat g_n$ its plug-in version, Step~2 and continuous mapping give
$\sup_u|\widehat g_n-g|\plim0$, and the decomposition of part (ii) of the
proof of Theorem~\ref{thm:mart}, applied \emph{mutatis mutandis} with this $g$, yields
$\widehat V^{0}_n\plim4\int_0^Tg\,d\Lam=V^{0}$.
\end{proof}

\phantomsection\label{prf:tstat}
\begin{proof}[Proof of Theorem~\ref{thm:tstat}]
All convergence arguments below are conducted under the regular conditional law given $\mathcal G$, for almost every volatility path. By Lemma~\ref{lem:driftremoval} the demeaned returns, and hence $\QCd_T$
and $\widehat V^{0}_n$, are pathwise identical to those computed from the
martingale part. We may therefore take $\alpha\equiv0$. Then
Lemma~\ref{lem:limitdm}, Theorem~\ref{thm:consistency}, and
Proposition~\ref{prop:condgaussdm} give $\QCd_T=L^{0}+o_p(1)$,
$\widehat V^{0}_n\plim V^{0}$, and $V^{0}>0$ a.s.. By Slutsky
$t^{0}_n=L^{0}/\sqrt{V^{0}}+o_p(1)$. Conditionally on $\Gsig$, $V^{0}$ is
$\Gsig$-measurable, so the ratio $L^{0}/\sqrt{V^{0}}$ is exactly standard
normal by Proposition~\ref{prop:condgaussdm}. Hence
$\Prob(t^{0}_n\le z\mid\Gsig)\to\Phi(z)$ for a.e.\ volatility path, and
unconditionally by dominated convergence.
\end{proof}

\phantomsection\label{prf:powerdm}
\begin{proof}[Proof of Corollary~\ref{cor:powerdm}]
By Lemma~\ref{lem:driftremoval} applied window by window, all demeaned
quantities are pathwise identical to those computed from the martingale
part, so we may take $\alpha\equiv0$ on $[0,KT]$. The proofs then go
through \emph{mutatis mutandis} with Proposition~\ref{prop:condgaussdm},
Lemma~\ref{lem:limitdm} and the demeaned part of Lemma~\ref{lem:Vhatalt}, and
Theorem~\ref{thm:tstat}
in place of Proposition~\ref{prop:condgaussraw}, Lemma~\ref{lem:limitraw}
and Theorem~\ref{thm:mart}. For integrability,
$|L^{0,(k)}|\le4T\overline\sigma^{2}\sup_{u\in\mathcal I_k}|M^{(k)}_u|$
since $|\Rt_u|\le2\sup|M^{(k)}|$. Also $h^{0}$ is continuous, strictly
decreasing with a unique zero $u^\ast\in(0,T)$
(Proposition~\ref{prop:condgaussdm}), so
$|h^{0}(u)|=|\Lam_{u^\ast}-\Lam_u|\ge\underline\sigma^{2}|u-u^\ast|$ and
$V^{0}\ge4\underline\sigma^{6}\int_0^T(u-u^\ast)^{2}du
\ge\tfrac13\underline\sigma^{6}T^{3}$ while the last step minimizing $\int_0^T(u-u^\ast)^2du=\tfrac13\{u^{\ast3}+(T-u^\ast)^3\}$
over $u^\ast\in(0,T)$ at $u^\ast=T/2$..
\end{proof}

\subsection{The consistency of the test statistic without \texorpdfstring{$H_0$}{H0}} \label{app: tstat consistency}
Under $H_0$, $t^{(k)}_n\plim\tau^{(k)}$ is Theorem~\ref{thm:mart}. Under
the alternative, $\QC_T\plim L$ still holds (Lemma~\ref{lem:limitraw}),
but step (ii) of the proof of Theorem~\ref{thm:mart} invokes
Corollary~\ref{cor:righttag} with the anticipative weight
$(\Lam_T-\Lam_u)^2$, which requires $H_0$. The following lemma removes
that requirement.
 
\begin{lemma}[Consistency of the variance estimators without $H_0$]
\label{lem:Vhatalt}
Under Assumptions~\ref{ass:sv}--\ref{ass:mart}, $\widehat V_n\plim V$ as
$n\to\infty$. Under Assumption~\ref{ass:sv} and a window-constant drift,
$\widehat c\plim c$ and $\widehat V^{0}_n\plim V^{0}$. In neither
statement is $H_0$ assumed. Consequently $t_n\plim L/\sqrt V$ and
$t^{0}_n\plim L^{0}/\sqrt{V^{0}}$ with or without $H_0$.
\end{lemma}
 
\begin{proof}
\emph{Step 1 (uniform consistency of cumulative squared returns).}
With $\widehat\Lam_n(t):=\sum_{t_j\le t}(\Delta R_j)^2$, the classical
convergence of realized quadratic variation for continuous semimartingales
\citep[Theorem~I.4.47]{jacod1987limit} gives $\widehat\Lam_n(t)\plim\Lam_t$
for each $t$ without $H_0$, and monotonicity plus continuity of $\Lam$
upgrade this to $\sup_{t\le T}|\widehat\Lam_n(t)-\Lam_t|\plim0$ exactly as
in step (i) of the proof of Theorem~\ref{thm:mart}. Step~1 of the proof of
Theorem~\ref{thm:consistency} is a pathwise bound that does not use
$H_0$, so also $\sup_{t\le T}|\widehat\Lam^{0}_n(t)-\Lam_t|\plim0$.
 
\emph{Step 2 (the estimators are Stieltjes integrals against
$\widehat\Lam_n$).} Regard $\widehat\Lam_n$ as the right-continuous step
function with jump $(\Delta R_j)^2$ at $t_j$. For continuous $\varphi$ on
$[0,T]$, the Lebesgue--Stieltjes integral is
$\int_0^T\varphi\,d\widehat\Lam_n=\sum_j\varphi(t_j)(\Delta R_j)^2$, so by
\eqref{eq:Vhatraw} and \eqref{eq:Vhatdm}, exactly,
\[
\widehat V_n=4\!\int_0^T\bigl(\widehat\Lam_n(T)-\widehat\Lam_n(u)\bigr)^2
d\widehat\Lam_n(u),\quad
\widehat c=\int_0^T\frac uT\,d\widehat\Lam^{0}_n(u),\quad
\widehat V^{0}_n=4\!\int_0^T\bigl\{\widehat\Lam^{0}_n(T)-\widehat\Lam^{0}_n(u)-\widehat c\bigr\}^2
d\widehat\Lam^{0}_n(u),
\]
while $V=4\int_0^T(\Lam_T-\Lam_u)^2d\Lam_u$, $c=\int_0^T(u/T)\,d\Lam_u$,
$V^{0}=4\int_0^T\{\Lam_T-\Lam_u-c\}^2d\Lam_u$.
 
\emph{Step 3 (a deterministic continuity property).} Let $F_n,F$ be
nondecreasing on $[0,T]$ with $F(0)=F_n(0)=0$, $F$ continuous,
$\sup_{t\le T}|F_n(t)-F(t)|\to0$, and let $\varphi_n\to\varphi$ uniformly
with $\varphi$ continuous. Then $\int_0^T\varphi_n\,dF_n\to\int_0^T\varphi\,dF$.
Indeed,
$\bigl|\int\varphi_n\,dF_n-\int\varphi\,dF_n\bigr|\le\sup|\varphi_n-\varphi|\,F_n(T)\to0$,
and $\int\varphi\,dF_n\to\int\varphi\,dF$ because uniform convergence of
$F_n$ to a continuous $F$ is weak convergence of the finite measures $dF_n$
to $dF$ on $[0,T]$ and $\varphi$ is bounded and continuous.
 
\emph{Step 4 (conclusion).} By the subsequence criterion for convergence in
probability, it suffices to argue along a subsequence on which Step~1
holds almost surely. On such a subsequence apply Step~3 with
$F_n=\widehat\Lam_n$, $F=\Lam$ and
$\varphi_n(u)=4(\widehat\Lam_n(T)-\widehat\Lam_n(u))^2\to\varphi(u)=4(\Lam_T-\Lam_u)^2$
uniformly (continuous mapping from Step~1); this gives
$\widehat V_n\to V$ a.s.\ along the subsequence, hence
$\widehat V_n\plim V$. For the demeaned case, first apply Step~3 with
$F_n=\widehat\Lam^{0}_n$ and $\varphi_n=\varphi=u/T$ to get $\widehat c\to c$,
then with $\varphi_n(u)=4\{\widehat\Lam^{0}_n(T)-\widehat\Lam^{0}_n(u)-\widehat c\}^2
\to4\{\Lam_T-\Lam_u-c\}^2$ uniformly, which gives $\widehat V^{0}_n\to V^{0}$.
 
Finally, $V\ge\tfrac43\underline\sigma^6T^3>0$ and
$V^{0}\ge\tfrac13\underline\sigma^6T^3>0$ by the bounds in
\eqref{eq:Vlower} and Corollary~\ref{cor:powerdm}, which use only
Assumption~\ref{ass:sv}, so Slutsky's lemma with
Lemmas~\ref{lem:limitraw} and~\ref{lem:limitdm} gives the last claim.
\end{proof}

\subsection{Derivation of the Heston Benchmark}\label{app:heston}

Throughout this appendix returns are the martingale part of
\eqref{XDGP} with $\lambda=0$, $dR_u=\sqrt{v_u}\,dW^x_u$, with
$dv_u=\kappa(\theta-v_u)\,du+\xi\sqrt{v_u}\,dW^v_u$,
$d\langle W^x,W^v\rangle_u=\rho\,du$, and $v_0=\theta$. With respect to the non-constant drift arising from the stochastic volatility, i.e., $\frac{1}{2}(\theta-v_u)$, the impact is
discussed in Appendix~\ref{app:voldrift}. Recall $\phi=\kappa T$. Under these dynamics all
moments of $v_u$ are finite and bounded on $[0,T]$ \citep{cox1985theory},
so the stochastic integrals appearing in the expectations below are true
martingales, in particular $\E\!\int_0^T R_u^2 v_u\,du<\infty$ and
$\E\!\int_0^T v_u^2\,du<\infty$, and their expectations vanish.

\subsubsection{The return--variance covariance function}

By the It\^o product rule,
\[
d(R_u v_u)\;=\;R_u\,dv_u+v_u\,dR_u+d\langle R,v\rangle_u ,
\qquad
d\langle R,v\rangle_u=\sqrt{v_u}\cdot\xi\sqrt{v_u}\,\rho\,du
=\rho\xi v_u\,du .
\]
Taking expectations, using that $R$ is a martingale
($\E[v_u\,dR_u]=0$, $\E[R_u]=0$) and that $v_0=\theta$ implies
$\E[v_u]=\theta$ for all $u$, the function $m(u)=\E[R_uv_u]$ solves the
linear ODE
\[
m'(u)\;=\;\kappa\theta\,\E[R_u]-\kappa\,m(u)+\rho\xi\,\E[v_u]
\;=\;-\kappa\,m(u)+\rho\xi\theta ,
\qquad m(0)=0 ,
\]
whose solution is \eqref{eq:mheston}.

\subsubsection{\texorpdfstring{The $\E[L]$ term}{The E[L] term}}

From \eqref{eq:meanQC} and Lemma~\ref{lem:limitraw},
$\E[L]=2\int_0^T m(u)\,du$. Since
$\int_0^T(1-e^{-\kappa u})\,du=T-(1-e^{-\phi})/\kappa$, equation
\eqref{eq:qcovraw} follows.

\subsubsection{\texorpdfstring{The $\E[L^{0}]$ term}{The E[L0] term}} \label{app:L0}

By Lemma~\ref{lem:driftremoval} the demeaned statistic under a
window-constant drift coincides pathwise with that of the martingale part,
so $\E[L^{0}]$ is computed on $dR_u=\sqrt{v_u}\,dW^x_u$. The bridge
reweighting of Section~\ref{sec:demean} is exact here: writing
$\Rt_u=R_u-(u/T)R_T$ and using the martingale property,
\[
\E\bigl[R_T\,v_u\bigr]
=\E\bigl[R_u v_u\bigr]+\E\bigl[(R_T-R_u)\,v_u\bigr]
=m(u)+\E\bigl[\E(R_T-R_u\mid\F_u)\,v_u\bigr]
=m(u),
\]
so that $\E[\Rt_u v_u]=(1-u/T)\,m(u)$ and
$\E[L^{0}]=2\int_0^T(1-u/T)\,m(u)\,du
=\frac{2\rho\xi\theta}{\kappa}\int_0^T(1-u/T)(1-e^{-\kappa u})\,du$.
Splitting $(1-e^{-\kappa u})$, the integral is the difference of the two
elementary integrals $\int_0^T(1-u/T)\,du=T/2$ and, 
\[
\int_0^T\Bigl(1-\frac{u}{T}\Bigr)e^{-\kappa u}\,du
=\frac{1}{\kappa}-\frac{1-e^{-\phi}}{\kappa^{2}T}.
\]
Their difference $\frac{T}{2}-\frac{1}{\kappa}+\frac{1-e^{-\phi}}{\kappa^2T}$,
multiplied by $\frac{2\rho\xi\theta}{\kappa}$, is \eqref{eq:qcovdm}.

\subsubsection{Positivity of the brackets} 
Both brackets are functions of
$\phi=\kappa T$. Use an auxiliary function
$h(z)=(1-e^{-z})/z=\int_0^1e^{-zt}\,dt\in(0,1)$ for $z>0$. The raw bracket
is $T\{1-h(\phi)\}>0$. The demeaned bracket equals $(1/\kappa)\,g(\phi)$
with $g(z)=z/2-1+h(z)$. Then $g(0)=0$ and
$g'(z)=\tfrac12-\int_0^1te^{-zt}\,dt>\tfrac12-\int_0^1t\,dt=0$ for $z>0$,
so $g(\phi)>0$ for all $\phi>0$. Hence both $\E[L]$ and $\E[L^{0}]$ vanish if and only if
$\rho=0$ (given $\xi,\theta>0$), and the limiting forms
$\rho\xi\theta T^2$, $\rho\xi\theta T^2/3$ (as $\phi\to0$, from
$g(z)=z^2/6+O(z^3)$) and $2\rho\xi\theta T/\kappa$,
$\rho\xi\theta T/\kappa$ (as $\phi\to\infty$) follow directly.

\subsubsection{Relation to the third cumulant} \label{app:3rdcumulant}

For the continuous martingale $R$ with $d[R]_u=v_u\,du$, It\^o's formula
gives $dR_u^3=3R_u^2\,dR_u+3R_uv_u\,du$, whence
\[
\kappa_3(R_T)=\E[R_T^3]=3\int_0^T m(u)\,du
=\tfrac{3}{2}\,\E\bigl[L\bigr]
=\frac{3\rho\xi\theta}{\kappa}\left(T-\frac{1-e^{-\phi}}{\kappa}\right),
\]
the Heston form of the \cite{NP2020} aggregation identity with respect to the 3rd cumulant. Note this is the 3rd cumulant of the martingale part only. The closed-form 3rd cumulant of the full return process additionally reflects the $(\mu-\frac{1}{2}v_u)$ drift and therefore differs, which can be found in \cite{ORS2022}.

\subsubsection{The stochastic-volatility drift: an exact pathwise
correction} \label{app:voldrift}
 
The simulated log-price \eqref{XDGP} carries the drift $\mu-\tfrac12v_u$.
Its constant part $\mu$ is removed exactly by demeaning
(Lemma~\ref{lem:driftremoval}). What remains is the It\^o term
$-\tfrac12v_u$, which is \emph{not} window-constant and therefore lies
outside the assumption under which Proposition~\ref{prop:condgaussdm} and
Theorem~\ref{thm:tstat} were derived. Although the calibration
$\mu=\theta/2$, $v_0=\theta$ sets the \emph{expected} drift to zero, the
realized drift $\tfrac12(\theta-v_u)$ is nonzero and time-varying within
every window. Rather than bounding its effect, we compute it exactly.
 
\paragraph{Pathwise decomposition.}
Write the Heston log-return as
\[
R_u=\mu u-\tfrac12\Lam_u+M_u,
\qquad
\Lam_u=\int_0^u v_s\,ds,\qquad M_u=\int_0^u\sqrt{v_s}\,dW^x_s .
\]
Bridging is linear, so with $\widetilde M_u=M_u-(u/T)M_T$,
\[
\Rt_u=R_u-\frac uT R_T
=\widetilde M_u-\tfrac12\Bigl(\Lam_u-\frac uT\Lam_T\Bigr).
\]
Inserting this into the infill limit of Lemma~\ref{lem:limitdm} and
using $\int_0^T\Lam_u\,d\Lam_u=\tfrac12\Lam_T^2$ together with the
definition of $c=\frac1T\int_0^Tu\,d\Lam_u$ in \eqref{eq:rep},
\begin{align}
L^{0}
&=2\int_0^T\Rt_u\,d\Lam_u
=2\int_0^T\widetilde M_u\,d\Lam_u
-\int_0^T\Bigl(\Lam_u-\frac uT\Lam_T\Bigr)d\Lam_u
\nonumber\\
&=\underbrace{2\int_0^T\widetilde M_u\,d\Lam_u}_{=:L^{0}_{M}}
\;-\;\tfrac12\Lam_T^{2}+\Lam_T\,c
\;=\;L^{0}_{M}+D,
\qquad
D:=\Lam_T\Bigl(c-\tfrac12\Lam_T\Bigr)
=\Lam_T\int_0^T\Bigl(\frac uT-\frac12\Bigr)v_u\,du .
\label{eq:Ddecomp}
\end{align}
The identity \eqref{eq:Ddecomp} is exact, path by path. Three
observations organize what follows.
\begin{enumerate}
\item[(i)] $L^{0}_{M}$ is the demeaned limit computed from the martingale
part alone, so Proposition~\ref{prop:condgaussdm} applies to it
\emph{mutatis mutandis}: $L^{0}_{M}=2\int_0^Th^{0}(u)\,dM_u$ and, under $H_0$,
$L^{0}_{M}\mid\Gsig\sim\mathcal N(0,V^{0})$ with $V^{0}$ as in
\eqref{eq:V}.
\item[(ii)] $D$ is a functional of $\Lam$ alone, hence
$\Gsig$-measurable. Conditionally on the volatility path the
stochastic-volatility drift is therefore a \emph{pure location shift}:
\begin{equation}\label{eq:L0shift}
L^{0}\mid\Gsig\;\sim\;\mathcal N\bigl(D,\;V^{0}\bigr)
\qquad\text{under }H_0 .
\end{equation}
Neither the conditional variance nor the conditional shape is affected.
\item[(iii)] $D$ is driven by the within-window \emph{tilt} of the
variance path. A rising variance path gives $D>0$, a falling one $D<0$, and constant
volatility gives $D=0$ exactly (then $c=\tfrac12\Lam_T$, cf.\
Remark~\ref{rem:bridge}).
\end{enumerate}
 
\paragraph{The closed form \eqref{eq:qcovdm} is exact at stationarity.}
By \eqref{eq:Ddecomp}, $\E[L^{0}]=\E[L^{0}_{M}]+\E[D]$, and $\E[L^{0}_{M}]$
is the quantity computed in Appendix~\ref{app:L0}. For $\E[D]$, write
$\E[v_uv_s]=\theta^2+\mathrm{c}(s,u)$ with the stationary CIR
autocovariance $\mathrm{c}(s,u)=\frac{\xi^2\theta}{2\kappa}e^{-\kappa|u-s|}$.
Then
\[
\E[D]=\int_0^T\!\!\int_0^T\Bigl(\frac uT-\frac12\Bigr)
\bigl\{\theta^2+\mathrm{c}(s,u)\bigr\}\,ds\,du
=\int_0^T\Bigl(\frac uT-\frac12\Bigr)
\underbrace{\frac{\xi^2\theta}{2\kappa^2}
\bigl(2-e^{-\kappa u}-e^{-\kappa(T-u)}\bigr)}_{=\int_0^T\mathrm{c}(s,u)\,ds}\,du
=0 ,
\]
because the $\theta^2$ term integrates to zero and the remaining
integrand is the product of a function antisymmetric about $u=T/2$ and
one symmetric about $u=T/2$. Hence $\E[L^{0}]=\E[L^{0}_{M}]$ and
\eqref{eq:qcovdm} holds without modification. 
Under the deterministic initialization $v_0=\theta$ used in the simulations, the covariance contains an additional transient term, and $\E[D]$ need not vanish. Table~\ref{tab:driftcorr} quantifies the resulting effect on the statistic and its rejection rates.
%
 
\paragraph{An exactly valid statistic for the Heston model.}
Because $D$ depends only on $\Lam$, it has a feasible counterpart built
from the same demeaned squared returns in
\eqref{eq:Vhatdm}:
\begin{equation}\label{eq:Dhat}
\widehat D
=\widehat\Lam^{0}_n(T)\Bigl(\widehat c-\tfrac12\widehat\Lam^{0}_n(T)\Bigr),
\qquad
t^{0,D}_n=\frac{\QCd_T-\widehat D}{\sqrt{\widehat V^{0}_n}} .
\end{equation}

\begin{proposition}[Exact correction under Heston dynamics]\label{prop:hestonD} 
Let returns follow \eqref{VDGP}--\eqref{XDGP} with $\lambda=0$ and
$\mu$ constant, and let $H_0$ hold. Then, as $n\to\infty$,
\[
t^{0,D}_n\;\dlim\;\mathcal N(0,1),
\qquad
t^{0}_n\;=\;t^{0,D}_n+\frac{\widehat D}{\sqrt{\widehat V^{0}_n}},
\qquad
\frac{\widehat D}{\sqrt{\widehat V^{0}_n}}\;\plim\;
\frac{D}{\sqrt{V^{0}}},
\]
the first convergence holding conditionally on $\Gsig$ for almost every
volatility path. In particular $t^{0}_n\mid\Gsig$ is asymptotically
$\mathcal N\bigl(D/\sqrt{V^{0}},\,1\bigr)$: the uncorrected statistic of
Theorem~\ref{thm:tstat} carries a path-specific location error
$D/\sqrt{V^{0}}$ and no scale error.
\end{proposition}
 
\begin{proof}
All convergence arguments below are conducted under the regular conditional law given $\mathcal G$, for almost every volatility path.
The drift $\mu u-\tfrac12\Lam_u$ is a continuous finite-variation path,
so $[\Rt]=[R]=\Lam$ and the convergences $\QCd_T\plim L^{0}$
(Lemma~\ref{lem:limitdm}), $\sup_t|\widehat\Lam^{0}_n(t)-\Lam_t|\plim0$,
$\widehat c\plim c$ and $\widehat V^{0}_n\plim V^{0}$ (proof of
Theorem~\ref{thm:consistency}) hold as stated: the window-constant-drift
assumption enters those results only through
Lemma~\ref{lem:driftremoval}, which is used to identify
$L^{0}$ with the martingale-part limit, not to establish the limits
themselves. Continuous mapping gives $\widehat D\plim D$, hence
$\QCd_T-\widehat D\plim L^{0}-D=L^{0}_{M}$ by \eqref{eq:Ddecomp}, and
$V^{0}>0$ a.s.\ by Proposition~\ref{prop:condgaussdm}. Slutsky then
yields $t^{0,D}_n=L^{0}_{M}/\sqrt{V^{0}}+o_p(1)$, which is exactly
standard normal conditionally on $\Gsig$ by observation~(i). The
decomposition of $t^{0}_n$ is immediate from \eqref{eq:Dhat}.
\end{proof}
 
\begin{remark}[Beyond Heston]\label{rem:itodrift}
The decomposition \eqref{eq:Ddecomp} uses only two facts: that the
log-price drift equals a constant minus $\tfrac12 v_u$, and that
$\int_0^T\Lam_u\,d\Lam_u=\tfrac12\Lam_T^2$. The first is not a Heston
feature but the It\^o correction of \emph{any} continuous price process
with constant proportional drift, $dS_u/S_u=\mu\,du+\sqrt{v_u}\,dW^x_u$.
Hence $t^{0,D}_n$ is asymptotically $\mathcal N(0,1)$ under $H_0$
throughout $\mathcal{SVSM}^c$ whenever the window-constant-drift
assumption is placed on the \emph{price} rather than on the log-price,
with no restriction on the variance dynamics. We nevertheless report
$t^{0}_n$ in the empirical work: the correction presumes a specific drift
convention that cannot be verified on the data, and the Heston
experiments show that its omission costs a location error of order
$10^{-2}$ on the $\mathcal N(0,1)$ scale.
\end{remark}

\end{document}